\documentclass[oneside,12pt]{amsart}
\usepackage{amssymb, amscd}
\usepackage[all]{xy}

\theoremstyle{plain}
\newtheorem{theorem}{Theorem}[section]
\newtheorem{lemma}[theorem]{Lemma}

\newtheorem*{theorem*}{}
\newtheorem{proposition}[theorem]{Proposition}
\newtheorem{corollary}[theorem]{Corollary}

\theoremstyle{definition}

\newtheorem{definition}[theorem]{Definition}

\theoremstyle{remark}
\newtheorem{remark}[theorem]{Remark}

\DeclareMathOperator{\C}{\mathbb C}
\DeclareMathOperator{\Z}{\mathbb Z}
\DeclareMathOperator{\R}{\mathbb R}
\DeclareMathOperator{\F}{\mathbb F}
\DeclareMathOperator{\soc}{soc}

\DeclareMathOperator{\Hom}{Hom}

\begin{document}
\title{The Cohomology of the Sylow $2$-Subgroup of the Higman-Sims Group}

\author{A. Adem}
\address{Mathematics Department\\
         University of Wisconsin\\
         Madison WI 53706} 
\email{adem@math.wisc.edu}
\thanks{The first author was partially supported by the NSF, NSA and
CRM-Barcelona}

\author{J. F. Carlson}
\address{Mathematics Department \\
         University of Georgia\\
         Athens GA 30602}
\email{jfc@math.uga.edu}
\thanks{The second author was partially supported by the NSF}

\author{D. B. Karagueuzian}
\address{Mathematics Department\\
         University of Wisconsin\\
         Madison WI 53706} 
\email{dikran@math.wisc.edu}
\thanks{The third author was 
partially supported by an NSF postdoctoral fellowship and the NSA}

\author{R. James Milgram}
\address{Mathematics Department \\
         Stanford University\\
	 Stanford CA 94305}
\email{milgram@math.stanford.edu}
\thanks{The fourth author was partially supported by the NSF and CRM-Barcelona}


\begin{abstract}
In this paper we compute the mod $2$ cohomology of the Sylow
$2$-subgroup of the Higman--Sims group $HS$, one of the $26$ sporadic
simple groups. We obtain its Poincar\'e series as well as an
explicit description of it
as a ring with $17$ generators and $79$ relations.
\end{abstract}

\maketitle
\section{Introduction}
Recently there has been substantial progress towards computing the mod $2$
cohomology of low rank sporadic simple groups. In fact the mod 2 cohomology
of every sporadic simple group not containing $(\mathbb Z/2)^5$ has been
computed, with the notable exceptions of $HS$ (the Higman--Sims group)
and $Co_3$ (one of the Conway groups). The reason for these exceptions
is that these two groups have large and complicated Sylow $2$-subgroups,
with many conjugacy classes of maximal elementary abelian subgroups.
The Higman--Sims group has order
$44,352,000=2^9\cdot 3^2\cdot 5^3\cdot 7\cdot 11$, and its largest elementary
abelian 2-subgroup is of rank equal to four. It is a subgroup of 
$Co_3$ of index equal to $11178=2\cdot 3^5\cdot 23$, hence $Syl_2(HS)$ is an
index $2$ subgroup of $Syl_2(Co_3)$.

In this paper we compute the mod $2$ cohomology ring of $S=Syl_2(HS)$,
obtaining explicit generators and relations. In a sequel we will
determine the necessary stability conditions for computing the
cohomology of $HS$ itself. This is a step towards obtaining
a calculation of the mod $2$ cohomology of 
$Co_3$. This latter group is of particular interest because of its
relation to a homotopy--theoretic construction due to Dwyer and Wilkerson
(see \cite{Be}) and it would seem that the most viable way of accessing
this group is via $HS$.

The calculation we present is long and highly technical, involving
techniques from topology, representation theory and computer algebra.
Our main result is the following

\begin{theorem}
\label{first-result}
The mod $2$ cohomology of $S=Syl_2(HS)$ has Poincar\'e series
\[
\frac{(1+x)^2(1-x+x^2)(1+2x-x^5)}{(1-x)^2(1-x^4)(1-x^8)}.
\]
As a ring, $H^*(S,\mathbb F_2$) has seventeen generators, in degrees
$$1,1,1,2,2,2,2,3,3,3,3,4,4,5,6,7,8,$$ 
and a (minimal) set of seventy--nine
relations.
\end{theorem}

\begin{remark} Nine of the generators above can be given as Stiefel-Whitney
classes associated to some of the irreducible representations of $S$
(see \S\ref{s:sw-class}), although
the remaining eight had to be determined by other means.  
The complete set of relations and the Steenrod operations on a representative
set of generators are
described in two appendices at the end of the paper.
\end{remark}

We briefly outline our method of proof. The group $S$ can be expressed
as a semidirect product $(\mathbb Z/4)^3:U_3$, where $U_3\cong D_8$
is the group of upper triangular $3\times 3$ matrices over $\mathbb F_2$.
Using a recent result in \cite{KS}, we calculate the $E_2$--term of
the Lyndon--Hochschild--Serre spectral sequence associated to this
extension. This involves using an explicit decomposition of the symmetric
algebra of a module. 

The second ingredient is a computer--assisted verification of the
cohomology of $S$ through degree ten. Combined with the structure
of the $E_2$ term we obtain the following crucial result

\begin{theorem}
The mod $2$ spectral sequence associated to $S=4^3:U_3$ collapses
at $E_2$, and yields the Poincar\'e series above.
\end{theorem}

A critical aspect of this is that the generators for the $E_2$
term occur in low degrees, where the computer can provide enough
information to ensure a collapse. This method seems to be
quite effective for many cohomology computations and gives a very
easy calculation of the mod $2$ cohomology of $U_4$.  Moreover,
with somewhat more effort it yields the mod $2$ cohomology of
$U_5$ as well.  This last result is particularly interesting since
$U_5$ is the Sylow $2$-subgroup of the two sporadics $M_{24}$
and $He$.

The next step consists of computing the restriction map from
$H^*(S,\mathbb F_2)$ to the cohomology of centralizers of
rank two elementary abelian subgroups, where we obtain an
image ring having the same Poincar\'e series as the cohomology.
 From this we conclude

\begin{theorem}
\label{second-result}
The cohomology of $S=Syl_2(HS)$ is detected by the centralizers of
rank two elementary abelian subgroups, and its explicit image
is the ring described above, with $17$ generators and $79$ relations.
\end{theorem}

\begin{remark}
There are nine such subgroups up to conjugacy in $S$.  
Explicit generators are described in \S\ref{s:sw-class} via their
restrictions to these nine centralizers.
\end{remark}

As mentioned above, the results here will be used in a sequel to
calculate the cohomology of $HS$, and from there obtain information
on the cohomology of $Co_3$, thus bringing us closer to a complete
understanding of rank $4$ sporadic simple groups. 

Throughout this paper coefficients will be in $\mathbb F_2$, so they
are suppressed. Occasionally $k$ will be used to denote this coefficient
ring, especially if it appears in a representation-theoretic context.
We refer the reader to \cite{AM} for background on group cohomology.

\section{The subgroup structure of $S = Syl_2(HS)$}

The Sylow $2$-subgroup of the Higman--Sims group (denoted $S$ from here on)
has a description as a semi--direct product
$S = 4^3\colon D_8$ with $v_1$, $v_2$, $v_3$
the generators of the $4^3$ while $D_8 = \{ t,s~|~t^2=s^4 = (ts)^2 = 1\}$, and
$$v_1^t~=~v_3^{-1}, ~~
v_2^t~=~v_2^{-1}, ~~
v_1^s ~=~v_2, ~~
v_2^s ~=~v_3, ~~
v_3^s ~=~v_2^{-1}v_1v_3.$$
In this section we will develop
the required details about the subgroup structure
of $S$ to enable us to
understand its cohomology. In our notation, if $x,y\in G$ are elements in a 
group, then $x^y=yxy^{-1}$.

To start we have

\begin{lemma}
$<v_1v_3> = \mathbb Z/4$ is a normal subgroup
of $S$ and $S/< v_1v_3 >\cong 2\wr 2\wr 2$.
\end{lemma}

\begin{proof}
In the quotient we have $v_1^t = v_3^{-1} \sim v_1$, while
$v_1^{ts^2} = v_1^{-1}$.  Also $v_2^t = v_2^{-1}$ while $v_2^{ts^2} = 
v_2v_3^{-1}v_1^{-1}$, which is equal to $v_2$ in the quotient. Thus
$$< v_1, ts^2> \cong D_8,~< v_2, t> \cong D_8$$
and these two copies of $D_8$ commute with each other, giving a copy of
$D_8\times D_8$ in the quotient.  Next, note that $ts$ exchanges $t$, $ts^2$,
and also exchanges $v_1$, $v_2^{-1}$, hence the two copies of $D_8$ above, and
the extension $(D_8\times D_8)\colon < ts> \cong 2\wr 2\wr 2$.
\end{proof}

\begin{definition}
The group $K_{\beta}$ is the inverse image in $S$
of the index two subgroup $D_8\times D_8\subset S/< v_1v_3>$.
\end{definition}

$K_{\beta}$ is given explicitly as an extension
$$< v_1, v_2, v_3>\colon
< t, s^2> = 4^3\colon 2^2.$$
Moreover, $S$ is the split extension $K_{\beta}\colon 2 = K_{\beta} \colon < ts
>$. 

Next we examine the maximal elementary abelian subgroups in $S$.
We show that there are precisely eight copies of
$2^4\subset S$.
To begin note that there are exactly 5 conjugacy classes of
maximal $2$-elementaries in
$2\wr 2\wr 2$, first $3$ conjugacy classes of $2^4$'s: 
$$2^4_I ~=~ < v_1^2, ts^2, v_2^2, t>,
~2^4_{II}~ =~<
v_1^2, v_1ts^2, v_2^2, v_2^{-1}t>$$
both normal, and $2^4_{I,II} = <
v_1^2 ts^2, v_2^2, v_2^{-1}t>$ with Weyl group $2^2$.  Then there
are two conjugacy classes of $2^3$'s, each with Weyl group $D_8$,
$$2^3_I~=~ < (v_1v_2)^2, s^2, ts>,
~2^3_{II}~ = ~< (v_1v_2)^2,
v_1v_2^{-1}s^2, ts>$$
which together generate a copy of $D_8\times 2$,
$$\Delta (D_8)\times 2 ~=~ < v_1v_2^{-1} s^2, ts>.$$
This is all standard and can be found in many references.  Lifting these
groups to $S$ we find

\begin{lemma}
The following list describes the groups above and their lifts:
$$2^4_I, ~~~D_8*D_8*4$$ 
$$2^4_{II}, ~~~D_8*D_8*4$$
$$2^4_{I,II}, ~~~(8\colon Aut(8))*4$$
$$\Delta (D_8)\times 2, ~~~D_8\times D_8$$
\end{lemma}

\begin{proof}
We note that the lift of $2^4_I$ is given as
$$< v_1v_3, v_1^2, ts^2, v_2^2, t> =
< v_1v_3^{-1}, ts^2, t, v_1^2, v_1v_3^{-1}v_2^2>$$
On the other hand, note that $v_1v_3^{-1}v_2^2$ commutes with $s^2$, $t$,
consequently with each of the remaining four generators, while the first two
generators commute with the third and fourth, and
$< v_1v_3^{-1}, ts^2> = D_8$, $< t, v_1^2> = D_8$,
and all three have the central element $(v_1v_3)^2$ in common.  The verification
for the second group is similar.  For the third, note that $(v_1s^2)^2 = v_1v_3$ so that
$v_1s^2$ has order $8$.  Also conjugation with $v_2^2$ takes the element to its fifth
power, while conjugation by $t$ takes it to its inverse.  This gives the extension
$8\colon Aut(8)$.  The final generator can again be choosen as $v_1v_3^{-1}v_2^2$.

It remains to check the lift of $\Delta (D_8)\times 2$.  Note, first that
$< s^2, v_1v_2^{-1}> = D_8$ as a subgroup of $S$, while
$v_1v_3$, $ts$ both commute with the elements of this $D_8$.  On
the other hand $< v_1v_3, ts> = D_8$ as well, and the final
statement follows.
\end{proof}

Now, to determine the structure of the set of $2^4$'s in $S$ we
check the inverse images of the conjugates of $2_I^3$ and $2_{II}^3$.  
The normalizer of each of these groups has order $2^6$ so there
are a total of $4$ such groups in $2\wr2\wr2$, with the remaining
two given by conjugation with $v_1$,
$$(2_I^3)^{v_1}~=~ 
< v_1v_2ts, v_2v_3^{-1}s^2, (v_1v_2)^2>,
~(2_{II}^3)^{v_1}~=~< v_1v_2ts, v_2v_3^{-1}s^2, (v_1v_2)^3>
$$
The lift of each to $S$ is a copy of $2^2\times D_8$ which contains
exactly two copies of $2^4$, so $S$ has
at least eight copies of $2^4$.  In fact we have

\begin{corollary}
There are no $2^5$'s contained in $S$.  
There are exactly eight copies of
$2^4$ contained in $S$, breaking up into three conjugacy classes,
two with two copies each and one with four.
\end{corollary}

\begin{proof}
The only thing that needs to be pointed out
is that if there were a $2^4$ which did not contain the central
element $(v_1v_3)^2$,  (and hence a $2^5$ obtained by adjoining $(v_1v_3)^2$),
then it would project non-trivially
to one of the four $2^4$'s in the quotient $2\wr 2\wr 2$.  But we
have identified the lifts of these groups as copies of groups having
rank three!  Hence, every possible $2^4$ contains $(v_1v_3)^2$,
and hence projects to one of the two conjugacy classes of extremal
$2^3$'s.  Consequently, it must lie in one or the other of the four
copies of $2^2\times D_8 \subset S$ that we have constructed.
\end{proof}

We now list the eight $2^4$'s explicitly.
$$2_A~=~< (v_1v_3)^2, (v_1v_2)^2,s^2, ts >$$
$$2_A^{v_1}~=~< (v_1v_3)^2, (v_1v_2)^2,v_1^{-1}v_3s^2, v_1v_2ts>$$
$$2_B~=~< (v_1v_3)^2, (v_1v_2)^2,v_1v_2^{-1}s^2,ts>$$
$$2_B^{v_1}~=~< (v_1v_3)^2, (v_1v_2)^2,v_2v_3^{-1}s^2,v_1v_2ts>$$
$$2_B^t~=~< (v_1v_3)^2, (v_1v_2)^2,v_3v_2^{-1}s^2,v_3v_2^{-1}ts>$$
$$2_B^{v_1t}~=~< (v_1v_3)^2, (v_1v_2)^2,v_1v_2^{-1}s^2,v_1v_3ts>$$
$$2_C~=~< (v_1v_3)^2, (v_1v_2)^2,s^2,v_1v_3ts>$$
$$2_C^{v_1}~=~< (v_1v_3)^2, (v_1v_2)^2,v_1v_3^{-1}s^2,v_2v_3^{-1}
ts>$$

\begin{lemma}
There is an outer automorphism $\alpha : S\to S$ which exchanges
$2_A$ and $2_C$.
\end{lemma}
\begin{proof}
We use the description of $S$ as $4^3\colon D_8$
where we write 
$$D_8 = 2^2\colon 2 = < ts, ts^{-1}> \colon
< t>.$$  
But we can replace this copy of $D_8$ by
$$< v_2v_3^{-1}ts, v_1v_2^{-1}ts^{-1}>\colon < t>$$
and $\alpha$ is defined as the identity on $4^3$ and the correspondence
above on $D_8$.
\end{proof}

We show that there are precisely two conjugacy classes of
$D_8\times D_8\subset S$.

\begin{remark}
The group $D_8\times D_8$ constructed in the proof
above as the lift of $\Delta (D_8)\times 2$ is
$$< s^2, v_1v_2^{-1}> \times < ts, v_1v_3>.$$
We can also construct a second copy of 
$$D_8\times D_8\subset S$$
as
$$< ts, v_2v_3^{-1}> \times< ts^{-1}, v_1v_2^{-1}>.$$
It is direct to check that the two $(D_8)^2$ above are not conjugate in
$S$.  In fact, the intersection of the second copy of $(D_8)^2$ and
$< v_1v_3 >$ is just $< (v_1v_3)^2>$ and its
image in $2\wr2\wr2$ is easily seen to be $D_8*D_8$.  Thus there are
at least two conjugacy classes of $D_8\times D_8$'s contained in $S$.  
Shortly we will show that there are exactly two.
\end{remark}

\begin{lemma}
The span $< 2_I^{e_1}, 2_J^{e_2}>$
is always one of the three groups $D_8\times D_8$, $2^2\times D_8$, or
$2^{2+4} = Syl_2(L_3(4))$.
\item{(a)}  It is $2^{2+4}$ if and only if $I = J$ and $e_2 = e_1v_1$.  
Consequently, there are exactly $3$ conjugacy classes of $2^{2+4}$'s
in $S$, two normal and one containing two elements.
\item{(b)}  It is $D_8\times D_8$ in case
$$< 2_A, 2_C^{v_1}> = < 2_B, 2_B^t>$$
and their conjugates by $v$, and also in the cases
$$< 2_A, 2_B^{v_1}> = < 2_C, 2_B^t>$$
and their four conjugates by $2^2 = < t, v_1>$.  
Consequently, there are exactly two conjugacy classes of $D_8\times D_8
\subset S$, one containing four groups and the other containing $2$.
\item{(c)}  In each of the remaining cases it is $2^2\times D_8$,
and each $2^2\times D_8$ is contained in a $D_8\times D_8$.
\end{lemma}

\begin{proof}
First we check the result for $2_A,~2_A^{v_1}$.  
We have
$$(s^2 v_1^{-1}v_3s^2)^2 = (v_1v_3)^2,~(s^2v_1v_2ts)^2=(v_2v_3)^2,$$
$$(tsv_1^{-1}v_3s^2)^2= (v_2v_3)^2,~(tsv_1v_2ts)^2 = (v_1v_2)^2,$$
and
this is a presentation of $2^{2+4}$.  The same calculations result for
$2_C$ using the automorphism above.  Moreover, $ts$ commutes with
$v_1v_2^{-1}$ while $s^2$ inverts it.  So this change cancels out
in the squares for
the pair $< 2_B, 2_B^{v_1}>$, and we have verified that
the groups asserted to be $2^{2+4}$'s in fact are.

The remaining statements are now easily checked by comparing with
the $D_8\times D_8$'s already constructed above.  But this gives
a complete list of possible pairs and the result follows.
\end{proof}

We show there there are exactly two conjugacy classes
of $2\wr2\wr2\subset S$.
Note that $t$ normalizes $< 2_A, 2_C^{v_1}>$, exchanging
the two copies of $D_8$, so
$$< 2_A, 2_C^{v_1}, t> \cong 2\wr 2\wr 2.$$
However, since $s^2$, $v_i^{e_1}v_j^{e_2}ts$ are not
conjugate in $S$, it follows that no $D_8\times D_8$ in the second
conjugacy class is contained in a $2\wr2\wr2$.

\begin{corollary}
There are exactly four copies of
$2\wr2\wr 2$ contained in $S$ forming two conjugacy classes
with conjugation by $v_1$ exchanging the groups in each class.
\end{corollary}

\begin{proof}
The normalizer of $< 2_B, 2_C^{v_1}>$
is obtained by adjoining $t$, $v_2^2$.  Thus there are three degree
two extensions of $< 2_B, 2_C^{v_1}>$ in $S$.  The extension
by $t$ is $2\wr 2\wr 2$.  Clearly, the extension by $v_2^2$ does
not give a $2\wr2\wr2$.

Finally, consider the extension by $tv_2^2$.  Replace the second
$D_8$ by 
$$< (v_2v_3)^2ts^{-1}, v_1v_2^{-1}>.$$  
Then
these two copies of $D_8$ commute with each other and their span
is $D_8\times D_8$.  Moreover, $tv_2^2$ exchanges them, giving an
isomorphism of this group with a second copy of $2\wr2\wr 2$.
\end{proof}

\begin{remark}
There is an automorphism of $S$ fixing
$< v_1, v_2, v_3>$ and exchanging the two conjugacy
classes of $2\wr2\wr2$'s constructed above.  Indeed, such an
automorphism can be given by setting
$$t ~\leftrightarrow ~ tv_2^2,
ts~\leftrightarrow ~ ts,
ts^{-1}~\leftrightarrow~ ts^{-1}(v_2v_3)^2$$
There are exactly six conjugacy classes of
maximal $2^3\subset S$.
\end{remark}

A computer analysis of $S$ using MAGMA or alternatively a direct analysis
of the two copies of $D_8*D_8*4$ and the $(8\colon Aut(8))*4$ obtained
as the lifts of the three $2^4$'s in $2\wr2\wr2$ shows that the maximal elementary
abelian subgroups of $S$ consist of the three $2^4$'s discussed above
and $6$ copies of $2^3$, all contained in the subgroup $K_{\beta}$:  one given
as $I = < v_1^2, v_2^2, v_3^2
>$ with centralizer $CI = 4^3$, and the other five all with
centralizer of the form $2^2\times 4$.  These centralizers are
$$CII~=~< v_1v_3^{-1}v_2^2, ts^2, v_3^2>$$
$$CIII~=~< v_1v_3^{-1}v_2^2, t, s^2>$$
$$CIV~=~< v_1v_3^{-1}v_2^2, s^2, v_1v_3t>$$
$$CV~=~< v_1v_3^{-1}v_2^2, v_1^2, v_1ts^2>$$
$$CVI~=~< v_1v_3^{-1}v_2^2, v_2t, v_2v_3^{-1}s^2>$$
with the $2^3$ subgroups denoted $II, \dots, VI$ respectively.

The group--theoretic information which we have described in this section will
be used subsequently to establish a detection theorem for
$H^*(S)$.

\section{Explicit Detection and Stiefel--Whitney Classes}
\label{s:sw-class}
In this section we describe the seventeen generators of $H^*(S)$
mentioned in \ref{first-result} in terms of their explicit
restrictions to the cohomology of the nine detecting subgroups, which
completely determines them in view of \ref{second-result}.  It turns out
that nine of the generators can be given as Stiefel-Whitney classes though
the remaining eight had to be determined by the computer.  For the nine 
Stiefel-Whitney classes we will be very explicit; for the remaining generators
we simply give convenient representatives. 

We use the 
following notation for the cohomology of the subgroups: 
$$H^*(2^m)=\mathbb F_2 [l_1,\dots ,l_m]$$
(where each $l_i$ is 1-dimensional);
$$H^*(4^3)=\mathbb F_2[b_1,b_2,b_3]
\otimes \Lambda (e_1,e_2,e_3)$$
(where the $e_i$ are 1-dimensional, and
$b_i$ is the Bockstein of $e_i$); and
$$H^*(2^2\times 4)=\Lambda (e)\otimes
\mathbb F_2 [l_2,l_3, b]$$
(where $|e|=1, b=\beta (e)$ and $|l_2|=|l_3|=1$).

One of the most powerful ways
of constructing elements in cohomology is to use
representations, i.e., homomorphisms $r_I\colon S \rightarrow
GL_n(\R )$, which, in turn, induce maps of classifying spaces,
$$B_{r_I}\colon B_S \rightarrow B_{GL_n(\R)}$$
and pull back the Stiefel-Whitney classes.  The restriction images
of these cohomology classes in the detecting groups are then
determined by taking the Stiefel-Whitney classes of the restrictions
of the $r_I$.  Here the determinations of the Stiefel-Whitney classes
are standard (see \cite{MM}, \cite{T}).

The real group-ring $\R(S)$ splits
into $33$ simple summands:
$$\R(S) = 8\R \oplus 12M_2(\R) \oplus M_2(\C) \oplus 8M_4(\R)
\oplus 3M_8(\R) \oplus M_8(\C),$$
but we will only need a small number of these representations.  

To begin consider the real representation $r_4$.  Here $r_4$ is
described by giving matrix images for the generators of the group:
$$v_1 \mapsto\left(\begin{array}{cccc}
-1&0&0&0\cr
0&1&0&0\cr
0&0&1&0\cr
0&0&0&-1\end{array}\right),~v_2\mapsto \left(\begin{array}{cccc}
1&0&0&0\cr
0&1&0&0\cr
0&0&-1&0\cr
0&0&0&-1\end{array}\right),~v_3 \mapsto \left(\begin{array}{cccc}
1&0&0&0\cr
0&-1&0&0\cr
0&0&-1&0\cr
0&0&0&1\end{array}\right)
$$
$$r_4(s) = \left(\begin{array}{cccc}
0&0&0&1\cr
1&0&0&0\cr
0&1&0&0\cr
0&0&1&0\end{array}\right),~ r_4(t) = \left(\begin{array}{cc}
K&0\cr
0&K\end{array}\right)
$$
where $K$ is the $2\times 2$ matrix $K = \left(\begin{array}{cc}
0&1\cr
1&0\end{array}\right)$. When
we restrict to $2_A$ we have that the first two generators map to
$I$, while $s^2 \mapsto \left(\begin{array}{cc}
0&I\cr
I&0\end{array}\right)$. Consequently,
when we diagonalize the image of $s^2$ to $\left(\begin{array}{cc}
-I&0\cr
0&I\end{array}\right)$, we see that $ts$ can also be diagonalized to
$$\left(\begin{array}{cccc}-1&0&0&0\cr
0&1&0&0\cr
0&0&1&0\cr
0&0&0&1\end{array}\right)$$
so the
representation becomes the sum of four one-dimensional
representations, and the total Stiefel-Whitney class becomes
$$(1 + l_3 + l_4)(1 + l_3) = 1 + l_4 + l_3(l_3+ l_4).$$
Next we consider $CI$.  Here we see at once that the total Stiefel-Whitney
class is
$$(1 + e_1)(1+e_3)(1+e_2+e_3)(1+e_1+e_2) = 1 + e_2(e_1+e_3).$$
As a final example of how to calculate these
restrictions, consider $CIII$.  Here, $v_1v_3^{-1}v_2^2 \mapsto
-I$ while the images of $s^2$ and $t$ have already been discussed.
Diagonalizing we have that the total Stiefel-Whitney class of the
restriction is
$$(1+e)(1+e+l_2)(1+e+l_3)(1+e+l_2+l_3) = 1 +
d_2(l_2,l_3) + (1+e)d_3(l_2, l_3).$$
Here, $d_2(x,y) = x^2 + xy +
y^2$ while $d_3(x,y) = x^2y + xy^2 = xy(x+y)$ are the Dickson
elements.

The other representations that will be needed are first, an
$8$-dimensional representation, $r_8$, which is very similar to the
$4$-dimensional representation considered above: $$v_1 \mapsto
\left(\begin{array}{cccc}J&0&0&0\cr
0&I&0&0\cr 0&0&I&0\cr
0&0&0&J\end{array}\right),~v_2 \mapsto
\left(\begin{array}{cccc}
I&0&0&0\cr
0&I&0&0\cr 0&0&J&0\cr
0&0&0&J\end{array}\right),~v_3 \mapsto
\left(\begin{array}{cccc}
I&0&0&0\\
0&J&0&0\\
0&0&J&0\\
0&0&0&I\end{array}\right)$$ 
$$s
\mapsto \left(\begin{array}{cccc}
0&0&0&I\cr
I&0&0&0\cr
0&I&0&0\cr
0&0&I&0\end{array}\right),~ t
\mapsto \left(\begin{array}{cccc}
0&K&0&0\cr
K&0&0&0\cr
0&0& 0&K\cr
0&0&K&0\end{array}\right)$$
where $I$ is the $2\times 2$ identity matrix and $J =
\left(\begin{array}{cc}
0&-1\cr
 1&0\end{array}\right)$. The remaining representations we need
are two-dimensional: the first of these, $r_{2,1}$, has the form
$$v_i \mapsto
J, ~1\le i \le 3, s \mapsto J,~t \mapsto K,
$$ the second, $r_{2,2}$, is $$v_i
\mapsto J, ~1\le i \le 3, s\mapsto I,~t \mapsto K,$$
and the third, $r_{2,3}$, is
$$v_1 \mapsto J, ~v_2\mapsto -J, v_3 \mapsto J,~s \mapsto K, t
\mapsto K.$$
Finally, we should mention the $1$-dimensional
representations which all factor in the form
$$S/S^{\prime} =
\langle v_1, s, t\rangle \cong 2^3 \longrightarrow \Z/2 = \{\pm 1\}$$
The first Stiefel-Whitney classes of $\langle v_1\rangle$,
$\langle s\rangle$, $\langle t\rangle$ respectively give the
$1$-dimensional generators for $H^*(S)$ while the second
Stiefel-Whitney classes of the three $2$-dimensional
representations above, together with the $4$-dimensional
representation, $r_4$,
give the two-dimensional generators.

The four $3$-dimensional generators do not occur as Stiefel-Whitney
classes, and were obtained by a computer calculation using MAGMA
as described previously.

The Stiefel-Whitney classes
$w_4$ and $w_8$ for the eight-dimensional representation
are generators. 

The $5$-dimensional generator is $Sq^2$
of one of the computer generated three-dimensional generators.  
The remaining $4$-dimensional generator, $n$, as well as the
seven-dimensional generator, $i$, also had to be determined by
MAGMA, and the $6$-dimensional generator can be given as
$Sq^2(n)$.

The following tables give the
restrictions of the generators discussed above.
\newpage
$$\mbox{\bf Table 1: Restrictions of One Dimensional Stiefel-Whitney Classes}$$
$$\bordermatrix{&w_1(v)&w_1(s)&w_1(t)\cr
CI&e_1+e_2+e_3&0&0\cr
CII&0&0&l_2\cr
CIII&0&0&l_2\cr
CIV&0&0&l_3\cr
CV&l_3&0&l_3\cr
CVI&l_2&0&l_2\cr
2_A&0&l_4&l_4\cr
2_B&0&l_4&l_4\cr
2_C&0&l_4&l_4 \cr}$$
$$\mbox{\bf Table 2: Restrictions of Two--Dimensional Stiefel-Whitney
Classes}$$
$$\bordermatrix{&w_2(r_{2,1})&w_2(r_{2,2})&w_2(r_{2,3})&w_2(r_4)\cr
CI&b_1+b_2+b_3&b_1+b_2+b_3&b_1+b_2+b_3&e_2(e_1+e_3)\cr
CII&el_2 +l_3(l_2+l_3)&el_2+l_3(l_2+l_3)&el_2 +
l_3(l_2+l_3)&l_2^2 \cr
CIII&el_2+l_3(l_2+l_3)&el_2+
l_3(l_2+l_3)&el_2&d_2(2,3) \cr
CIV&el_3 +
l_2(l_2+l_3)&el_3+l_2(l_2+l_3)&el_3&d_2(2,3)\cr
CV&el_3+l_2(l_2+l_3)&el_3+l_2(l_2+l_3)&el_3+l_2(l_2+l_3)&l_3^2\cr
CVI&el_2+l_3(l_2+l_3)&el_2+l_3(l_2+l_3)&el_2+l_3(l_2+l_3)&d_2(2,3)\cr
2_A&l_3(l_3+l_4)&0&0&l_3(l_3+l_4)\cr
2_B&l_3(l_3+l_4)&0&l_3^2&l_3(l_3+l_4) \cr
2_C&l_3(l_3+l_4)&0&l_4^2&l_3(l_3+l_4)\cr}$$
$$\hbox{\bf Table 3: Restrictions of $w_4$ and $w_8$ for the Eight-Dimensional Representation.}$$
$$\bordermatrix{&w_4&w_8\cr
CI&(b_1+b_2+b_3)^2 + b_2(b_1+b_3)&b_1b_3(b_1+b_2)(b_2+b_3)\cr
CII&d_2(2,3)^2&L\cr
CIII&d_2(2,3)^2&L\cr
CIV&d_2(2,3)^2&L\cr
CV&d_2(2,3)^2&L\cr
CVI&d_2(2,3)^2&L\cr
2_A&d_4(2,3,4)&M\cr
2_B&d_4(2,3,4)&M\cr
2_C&d_4(2,3,4)&M \cr}$$
where $d_2(2,3) = l_2^2 + l_2l_3 + l_3^2$,
$d_3(2,3) = l_2l_3(l_2+l+3)$, $d_4(2,3,4)$ is the fourth
Dickson invariant in the three classes $l_2$, $l_3$, and $l_4$, while
\begin{eqnarray*}
L &=&b^4 +b^2 d_2(2,3)^2 + bd_3(2,3)^2\\
M &=&l_1^8+l_1^4d_4(2,3,4) + l_1^2Sq^2(d_4(2,3,4)) +
l_1Sq^3(d_4(2,3,4)).
\end{eqnarray*}
Next we give the restrictions of the
computer generated indecomposables, $s$, $r$, $p$, and $q$,
in dimension three.
$$\bordermatrix{&s&r\cr
CI&e_1 e_2 e_3& e_1 e_2 e_3\cr
CII&0& e_1 l_2^2\cr
CIII&0&  e_1 l_2^2 + e_1 l_2 l_3 + e_1 l_3^2 + l_2^2 l_3 + l_2 l_3^2\cr
CIV&0& e_1 l_2^2 + e_1 l_2 l_3 + e_1 l_3^2\cr
CV&e_1 l_3^2&0\cr
CVI&e_1 l_2^2 + e_1 l_2 l_3 + e_1 l_3^2 + l_2 l_3^2 + l_3^3&0\cr
2_A&0&0\cr
2_B&l_3^3 + l_3^2 l_4&0\cr
2_C&0& l_3^2 l_4 + l_3 l_4^2\cr}
$$
$$\bordermatrix{&q&p\cr
CI& e_1 e_2 e_3 + e_2 b_3 + e_3 b_2& e_1 e_2 e_3\cr
CII&0& e_1 l_2^2\cr
CIII&0&e_1 l_2^2 + l_2^2 l_3 + l_2 l_3^2\cr
CIV&0& e_1 l_3^2\cr
CV& e_1 l_3^2&0\cr
CVI&e_1 l_2^2 + e_1 l_2 l_3 + e_1 l_3^2 + l_2 l_3^2 + l_3^3
    & e_1 l_2 l_3 + e_1 l_3^2\cr
2_A&0&
    l_1^2 l_4 + l_1 l_4^2 + l_2^2 l_3 + l_2 l_3^2 + l_3^2 l_4 + l_3 l_4^2\cr
2_B& l_3^3 + l_3^2 l_4&
    l_1^2 l_4 + l_1 l_4^2 + l_2^2 l_3 + l_2 l_3^2 + l_3^2 l_4 + l_3 l_4^2\cr
2_C&0& l_1^2 l_4 + l_1 l_4^2 + l_2^2 l_3 + l_2 l_3^2\cr}
$$
Note here that the computer has not always picked the simplest choices.  
For example $s + q$ restricts to $e_2b_3 + e_3b_2$ in $CI$ and $0$ in
the remaining $8$ centralizers.

Here is the expansion of the computer generated indecomposable, $n$
in dimension four:
$$\bordermatrix{&n\cr  
CI& e_1 e_2 b_3 + e_1 e_3 b_2 + e_2 e_3 b_1 + e_2 e_3 b_3\cr
CII&e_1 l_2^3 + e_1 l_2^2 l_3 + e_1 l_2 l_3^2\cr
CIII&e_1 l_2^3 + l_2^3 l_3 + l_2^2 l_3^2\cr
CIV&e_1 l_2^2 l_3 + e_1 l_2 l_3^2 + e_1 l_3^3\cr
CV&e_1 l_2^2 l_3 + e_1 l_2 l_3^2\cr
CVI&e_1 l_2 l_3^2 + e_1 l_3^3 + l_2^2 l_3^2 + l_3^4\cr
2_A&0\cr
2_B&l_1^2 l_3 l_4 + l_1 l_3 l_4^2 + l_2^2 l_3^2 + l_2 l_3^3
+ l_3^4 + l_3^3 l_4 + l_3^2 l_4^2\cr
2_C&l_1^2 l_4^2 + l_1 l_4^3 + l_2^2 l_3 l_4 + l_2 l_3^2 l_4 + l_4^4\cr}$$
This element, together with $w_4$ for the eight-dimensional representation
can be taken as the indecomposable generators in dimension $4$.
Finally, for the generator $i$ in degree $7$ the
restrictions are given as follows. For $CI$ the restriction is given as
\medskip

$\phantom \qquad e_1 e_2 e_3 b_3^2 + e_1 b_2^2 b_3 + e_1 b_2 b_3^2
+ e_2 b_1^2 b_3 + e_2 b_1 b_3^2 +  e_2 b_2^2 b_3
 + e_2 b_2 b_3^2 + e_3 b_1^2 b_2 + e_3 b_1 b_2^2 + e_3 b_2^3 + e_3 b_2^2 b_3
        + e_3 b_3^3$
\medskip

\noindent while for $CII$ through $CVI$ we have
\medskip

$$\bordermatrix{&i\cr
CII&    e_1 l_2^6 + e_1 l_2^5 l_3 + e_1 l_2^4 l_3^2 + l_2^4 l_3^3 + l_2^3 l_3^4 + l_2^2 l_3^5 + 
        l_2 l_3^6\cr
CIII&    e_1 l_2^6 + e_1 l_2^4 l_3^2 + e_1 l_2^3 l_3^3 + e_1 l_2 l_3^5 + e_1 l_3^6 + l_2^6 l_3 + 
        l_2^4 l_3^3 + l_2^2 l_3^5 + l_2 l_3^6\cr
CIV&    e_1 l_2^6 + e_1 l_2^5 l_3 + e_1 l_2^3 l_3^3 + e_1 l_2 l_3^5 + e_1 l_3^6\cr
CV&    e_1 l_2^2 l_3^4 + e_1 l_2 l_3^5 + e_1 l_3^6\cr
CVI&    e_1 l_2^6 + e_1 l_2^2 l_3^4 + e_1 l_2 l_3^5 + l_2^2 l_3^5 + l_2 l_3^6\cr}
$$
\medskip

\noindent   The restrictions to $2_A$, $2_B$, and $2_C$ are very 
long and complicated.  
For $2_A$ we obtain
\medskip

$\phantom \qquad
l_1^2 l_2^4 l_4 + l_1^2 l_2^2 l_3^2 l_4 + l_1^2 l_2^2 l_3 l_4^2 + l_1^2 l_2^2 l_4^3 + 
        l_1^2 l_2 l_3^2 l_4^2 + l_1^2 l_2 l_3 l_4^3 + l_1 l_2^4 l_4^2
+ l_1 l_2^2 l_3^2 l_4
 +
        l_1 l_2^2 l_3 l_4^3 + l_1 l_2^2 l_4^4
 + l_1 l_2 l_3^2 l_4^3 + l_1 l_2 l_3 l_4^4 + 
        l_2^6 l_3 + l_2^5 l_3^2 + l_2^4 l_3^3 + l_2^4 l_3^2 l_4 + l_2^4 l_3 l_4^2 + 
        l_2^4 l_4^3 + l_2^3 l_3^4 + l_2^2 l_3^4 l_4 + l_2^2 l_3^3 l_4^2 + l_2^2 l_3^2 l_4^3 +
        l_2^2 l_3 l_4^4 + l_2^2 l_4^5 + l_2 l_3^2 l_4^4 + l_2 l_3 l_4^5$
\medskip

\noindent
For $2_B$ the restriction is

\medskip

$\phantom\qquad
l_1^4 l_3^2 l_4 + l_1^2 l_2^4 l_4 + l_1^2 l_2^2 l_3^2 l_4 + l_1^2 l_2^2 l_3 l_4^2 +  l_1^2 l_2^2 l_4^3 + l_1^2 l_2 l_3^2 l_4^2 + l_1^2 l_2 l_3 l_4^3 + l_1^2 l_3^4 l_4 +  l_1^2 l_3^2 l_4^3 + l_1 l_2^4 l_4^2 + l_1 l_2^2 l_3^2 l_4^2 + l_1 l_2^2 l_3 l_4^3 + l_1 l_2^2 l_4^4 + l_1 l_2 l_3^2 l_4^3 + l_1 l_2 l_3 l_4^4 + l_1 l_3^4 l_4^2 + l_2^6 l_3  + l_2^5 l_3^2 + l_2^4 l_3^3 + l_2^4 l_4^3 + l_2^3 l_3^4 + l_2^2 l_3^4 l_4 + l_2^2 l_3^3 l_4^2 + l_2^2 l_3^2 l_4^3 + l_2^2 l_4^5 + l_2 l_3^5 l_4 + l_2 l_3^4 l_4^2  + l_2 l_3 l_4^5 + l_3^6 l_4 + l_3^3 l_4^4$
\medskip

\noindent
For $2_C$ the restriction is
\medskip

$\phantom\qquad
l_1^4 l_4^3 + l_1^2 l_2^4 l_4 + l_1^2 l_2^2 l_3^2 l_4 + l_1^2 l_2^2 l_3 l_4^2 + 
        l_1^2 l_2^2 l_4^3 + l_1^2 l_2 l_3^2 l_4^2 + l_1^2 l_2 l_3 l_4^3 + l_1^2 l_3^2 l_4^3 + 
        l_1^2 l_3 l_4^4 + l_1 l_2^4 l_4^2 + l_1 l_2^2 l_3^2 l_4^2 + l_1 l_2^2 l_3 l_4^3 + 
        l_1 l_2^2 l_4^4 + l_1 l_2 l_3^2 l_4^3 + l_1 l_2 l_3 l_4^4 + l_1 l_3^2 l_4^4 + 
        l_1 l_3 l_4^5 + l_1 l_4^6 + l_2^6 l_3 + l_2^5 l_3^2 + l_2^4 l_3^3 + l_2^4 l_3 l_4^2 + 
        l_2^3 l_3^4 + l_2^2 l_3^2 l_4^3 + l_2^2 l_3 l_4^4 + l_2 l_3^4 l_4^2 + l_2 l_3^3 l_4^3 
        + l_2 l_3^2 l_4^4 + l_3^6 l_4 + l_3^5 l_4^2 + l_3^4 l_4^3 + l_3^3 l_4^4$
\medskip

\noindent In the next section we begin the proof of our main theorem.

\section{Preliminaries on Modules}
\label{s:dec-SM}
We now turn to the proof of \ref{first-result}.  As our intial step
we have to determine $H^*(4^3)$ as a module over $\F_2(U_3)$.

Let $U_3\cong D_8$ denote the Sylow $2$-subgroup of $L_3(2)$.
In this section we describe the symmetric algebra of the natural $U_3$
module $M$ by giving a ``factorization'' of this algebra. 
This information will be used to describe the $E_2$ term of a spectral
sequence converging to $H^*(S)$. In order to
describe this factorization we must make a number of definitions and
recall a few results. We will be using ideas and methods from
\cite{KS}. We write the symmetric algebra as
$k[x,y,z]$, where the action of $U_3$ on the homogeneous polynomials
of degree $1$ preserves the flag of subspaces $< x >
\subset < x, y > \subset < x, y, z >$. 
Let's denote by $b,c,d$  the elements of $U_3$ which send $y \mapsto
y+x$, $z \mapsto z+x$, and $z \mapsto z+y$, respectively. ($b$ is
supposed to fix $z$, $c$ to fix $y$ and $d$ to fix $x$.) We write
$H_1$ for the subgroup $< b,c >$ and $H_2$ for the
subgroup $< d >$.

\begin{remark}
It is worth noting that the module $M$ is \emph{not} isomorphic to
$M^*$, but that the two differ by an outer automorphism of the group
$U_3$.  This means that the symmetric algebras $S^*(M)$ and $S^*(M^*)$
also differ by an outer automorphism of the group; thus if we are only
interested in Poincar\'e series it doesn't matter which one we study. 
However, it is \emph{not} true that $S^*(M)^* \approx S^*(M^*)$.
\end{remark}

\begin{definition} $d_1 := x, \qquad d_2 := y^2+xy, \qquad d_4 :=
z(z+y)(z+x)(z+y+x)$ 
\end{definition}

The following result is well-known:
 
\begin{proposition}
$k[x,y,z]^{U_3} = k[d_1,d_2,d_4]$
\end{proposition}

\begin{definition}$ \phi := z^2+yz, \qquad \theta := z^2+xz$. 
We denote by $N$, $K$,  the $U_3$-submodule of the symmetric algebra
generated by $\phi$, $\theta$, respectively.
\end{definition}

\begin{lemma}
$N$ is isomorphic to the permutation module $k[U_3/H_2]$ and has socle
$d_1^2$. $K$ is isomorphic to the permutation module $k[U_3/H_1]$ and
has socle $d_2$.
\end{lemma}

We will now exhibit three $U_3$-submodules of $k[x,y,z]$.  
\begin{definition}
\[
A := k[d_4], \qquad B := k \oplus K \oplus d_2K \oplus d_2^2K \oplus
\cdots 
\]
\[
C := k \oplus M \oplus N \oplus d_1N \oplus d_1^2N \oplus d_1^3N
\oplus \cdots
\]
\end{definition}

Note that $A$ is a $k[d_4][U_3]$-submodule of $k[x,y,z]$, while $B$ is
a $k[d_2][U_3]$-submodule and $C$ a $k[d_1][U_3]$-submodule. (The only
point to check is that $d_1M \subset N$, but this can be done by
working directly with the generating polynomials.)

Now we can describe the ``factorization'' mentioned above:

\begin{proposition}
There is an isomorphism of $U_3$-modules
\[
A\otimes B\otimes C \rightarrow k[x,y,z]
\]
induced by multiplication in the symmetric algebra.
\end{proposition}

\begin{proof}
For any triple of indecomposable modules $A_i$, $B_j$, $C_l$ appearing
in the direct sum decompositions of $A$, $B$, and $C$ above, we have
$\soc(A_i \otimes B_j \otimes C_l) = \soc(A_i) \otimes \soc(B_j)
\otimes \soc(C_l)$.  Thus $\soc(A \otimes B \otimes C) = \soc(A)
\otimes \soc(B) \otimes \soc(C) = k[d_4] \otimes k[d_2] \otimes
k[d_1]$. This shows that the map above is an isomorphism on the socle,
hence it is injective.  The Poincar\'e series of the two sides are
equal, by a calculation, so we have an isomorphism.
\end{proof}

\begin{remark}
No analogous factorization exists for larger fields or for more
variables.
\end{remark}

We can rewrite this factorization as a direct sum decomposition

\begin{proposition}
\label{p:dec-SM}
There is an isomorphism of $U_3$-modules:
\begin{align*}
k[x,y,z] & = & &  &  & k[d_4] \otimes (k \oplus M)  \\
         &   & &\oplus   &  & k[d_1,d_4] \otimes N  \\
         &   & &\oplus   &  & k[d_2,d_4] \otimes (K \oplus W)  \\
         &   & &\oplus   &  & k[d_1,d_2,d_4] \otimes F.
\end{align*}
Here $F = N \otimes K$ is a free $U_3$-module, appearing in degree
$4$, while $k$, $M$, $N$, $W = M \otimes K$, and $K$ appear in degrees
0,1,2,3, and 2 respectively.
\end{proposition}

\section{The Sylow $2$-subgroup of the
Higman-Sims group: Modules and Cohomology}

As we mentioned, the Sylow $2$-subgroup can be written as the semidirect product 
$({\mathbb Z_4})^3 \rtimes U_3$, where the action of $U_3$ on the
vector space $H^1(({\mathbb Z_4})^3)$ gives this cohomology group the
structure of the natural $U_3$-module $M$. It follows that as a
$U_3$-module, $H^*(({\mathbb Z_4})^3) = \Lambda^*(M) \otimes S^*(M)$. 

In this section we compute the Poincar\'e series of $H^*(U_3;
\Lambda^*(M) \otimes S^*(M))$.  We will also try to find a bound on
the degrees of generating elements.  To do this we must compute the
Poincar\'e series of the cohomologies of all the $U_3$-modules
appearing in the decomposition of $\Lambda^*(M) \otimes S^*(M)$. 
Since $\Lambda^*(M) = k \oplus M \oplus M^* \oplus k$, it is enough to
compute the cohomologies of the modules appearing in the decomposition
(\ref{p:dec-SM}) of $S^*(M)$, plus the cohomologies of the tensor
products of these modules with $M$ and $M^*$.  

The first thing to understand is how the tensor products of the
modules $k$, $M$, $N$, $K$, $W$, and $F$ with the modules $M$ and
$M^*$ break up into direct sums of indecomposables.

The answers are displayed in the table below.  A couple of new modules
appear; they are described below the table.

\bigskip
\begin{center}
\begin{tabular}{c | c | c }
$\otimes$ & $M$ & $M^*$ \\ \hline
 $k$ & $M$ & $M^*$  \\
 $M$ & $Y_9$ & $F \oplus k$ \\
 $N$ & $F \oplus N$ & $F \oplus N$\\
 $K$ & $W$ & $W^*$\\
 $W$ & $F \oplus X_{10}$ & $F^{\oplus 2} \oplus K$ \\
 $F$ & $F^{\oplus 3}$ & $F^{\oplus 3}$ \\
\end{tabular}
\end{center}
\bigskip 

 The module $X_{10}$ fits in an exact sequence $X_{10} \hookrightarrow
F \oplus F \twoheadrightarrow W$.  The module $Y_9$ is by definition
$M\otimes M$, and there is an exact sequence $M \hookrightarrow T
\oplus F \twoheadrightarrow Y_9$. This second exact sequence is
constructed by taking the exact sequence $k \hookrightarrow T
 \twoheadrightarrow M$ and applying $- \otimes_k M$.  Here we note
that $T$ is the permutation 
module $U_3 / < b >$.

 The next step is to compute cohomologies for each of the modules
listed in the table; we will need Poincar\'e series and degrees of
generators for the cohomologies $H^*(U_3, X)$ as modules over the
cohomology $H^*(U_3,k)$, for each indecomposable module $X$. 
We display the results in the form of a table and give the methods of
proof afterwards. 

\bigskip
\begin{center}
\begin{tabular}{c | c | c }
Module & Poincar\'e Series & Bound on Generator Degrees \\ \hline
 $k$ & $\frac{1}{(1-x)^2}$ & 0 \\
$M$ & $\frac{1}{(1-x)^2}$  & 1 \\
 $M^*$ & $\frac{1}{(1-x)^2}$  & 1 \\
 $N$ &  $\frac{1}{(1-x)}$  & 0 \\
 $K$ & $\frac{1}{(1-x)^2}$  & 1 \\
 $W$ & $1+ \frac{x}{(1-x)^2}$ & 2 \\
 $W^*$ & $ \frac{2-x}{(1-x)^2}$ & 1 \\
 $F$ & 1 & 0 \\
 $Y_9$ & $\frac{2-2x-x^2}{(1-x)^2}$ & 2 \\
 $X_{10}$ & $2 + x  + \frac{x^2}{(1-x)^2}$  & 3 \\
\end{tabular}
\end{center}
\bigskip 

Now we present brief arguments for these cohomology computations. 

In the case of a permutation module, i.e. the cases $k$, $N$, $K$, and
$F$, the cohomology is just the cohomology of the appropriate
subgroup, regarded as a module over the cohomology of $U_3$ via the
restriction map.  The degrees of the module generators can be worked
out from a knowledge of the restriction image. 

Now we turn to the more subtle cases.  First let us note that there is
an exact sequence $k \hookrightarrow T \twoheadrightarrow M$.  In
fact, more is true, there is also an exact sequence $M \hookrightarrow
N \twoheadrightarrow k$, and by dualizing we can get exact sequences
for $M^*$.  We'll just work with the first exact sequence. 

\begin{lemma}
The cohomology of $M$, i.e. $H^*(U_3,M)$, has Poincar\'e series
$(1-x)^{-2}$ and is generated by elements in degrees $0$ and $1$. 
\end{lemma}
\begin{proof}
Consider the long exact sequence in cohomology arising from the short
exact sequence $k \hookrightarrow T \twoheadrightarrow M$. The maps
$H^i(U_3,k) \rightarrow H^i(U_3,T)$ are just the restriction
homomorphisms for $< b > = {\mathbb Z}_2 \subset U_3$,
which can be shown to be surjective in all degrees $i$. This implies
that $H^i(U_3,M)$ is isomorphic to the  kernel of the restriction map
$H^{i+1}(U_3,k) \rightarrow H^{i+1}(U_3,T)$ for all $i \geq 0$. This
means that we can write the Poincar\'e series of the cohomology of $M$
as 
\[
x^{-1}\biggl( \frac{1}{(1-x)^2} - \frac{1}{(1-x)} \biggr) = \frac{1}{(1-x)^2}.
\]
The module generators may be taken to be in degrees $0$, $1$. 
\end{proof}

\begin{lemma}
The cohomology of $M^*$, i.e. $H^*(U_3,M^*)$, has Poincar\'e series
$(1-x)^{-2}$ and is generated by elements in degrees $0$ and $1$. 
\end{lemma}
\begin{proof}
$M$ and $M^*$ are interchanged by an outer automorphism of $U_3$. 
\end{proof}

\begin{lemma}
$H^*(U_3,W)$, has Poincar\'e series
$1+x(1-x)^{-2}$ and is generated by elements in degrees $\leq 2$. 
\end{lemma}
\begin{proof}
There is an exact sequence $W \hookrightarrow F \twoheadrightarrow K$.
\end{proof}

\begin{lemma}
$H^*(U_3,W^*)$, has Poincar\'e series
$(2-x)(1-x)^{-2}$ and is generated by elements in degrees $\leq 1$. 
\end{lemma}
\begin{proof}
There is an exact sequence $K \hookrightarrow F \twoheadrightarrow W^*$.
\end{proof}

\begin{lemma}
$H^*(U_3,Y_9)$, has Poincar\'e series
$(2-2x+x^2)(1-x)^{-2}$ and is generated by elements in degrees $\leq 2$. 
\end{lemma}
\begin{proof}
There is an exact sequence $M \hookrightarrow T\otimes M^*
\twoheadrightarrow Y_9$, and we consider the associated long exact
sequence in cohomology.  Since $T\otimes M^* \approx T \oplus F$, we
see that $H^i(U_3,T\oplus F) = H^i(U_3,T)$ is one-dimensional if $i >
0$. Furthermore, the 
induced maps $H^i(U_3,M) \rightarrow H^i(U_3,T)$ are surjective for
all $i > 0$.  This, plus the fact that the socle of $Y_9$ is
two-dimensional, determines the Poincar\'e series.  For the
information on the degrees of the generators, we must study the kernel
of the map $H^*(U_3,M) \rightarrow H^*(U_3,T)$. 
This kernel is generated by elements in degrees $\leq
2$.
\end{proof}

\begin{lemma}
$H^*(U_3,X_{10})$, has Poincar\'e series
$2+x+x^2(1-x)^{-2}$ and is generated by elements in degrees $\leq 3$. 
\end{lemma}
\begin{proof}
Use the exact sequence $X_{10} \hookrightarrow
F \oplus F \twoheadrightarrow W$. 
\end{proof}

\section{The $E_2$-term for the Sylow $2$-subgroup of the
Higman-Sims Group}

The algebraic computations in the preceding section will now be assembled
to describe the $E_2$--term of 
the Lyndon-Hochschild-Serre spectral sequence associated to the group
extension $S=4^3:D_8$. It can of course be described precisely as
$H^*(U_3,\Lambda^*(M)\otimes S^*(M))$.

To start we make a table of the modules appearing in the $E_2$-term, which
we regard as $\Lambda^*(M) \otimes S^*(M)$, and the Poincar\'e series
of their cohomologies. 

\bigskip

\begin{center}
\begin{tabular}{c | c | c | c}
Module & Decomposition & Poincar\'e Series & Degree Bound \\ \hline
$k$  & $k$ & $\frac{1}{(1-x)^2}$ & 0\\
$M$  & $M$ & $\frac{1}{(1-x)^2}$ & 1 \\
$M^*$  & $M^*$ & $\frac{1}{(1-x)^2}$ & 1\\
$M \otimes M$   &  $Y_9$ & $\frac{2-2x+x^2}{(1-x)^2}$ & 2\\
$M \otimes M^*$ &  $F \oplus k$ & $1 + \frac{1}{(1-x)^2}$ & 0\\
$N$ & $N$ & $\frac{1}{(1-x)}$ & 0\\
$N \otimes M$ &  $F \oplus N$ & $1 + \frac{1}{(1-x)}$ & 0 \\
$N \otimes M^*$ & $F \oplus N$ & $1 + \frac{1}{(1-x)}$ & 0 \\
$K$ & $K$ & $\frac{1}{(1-x)^2}$ & 1 \\
$K \otimes M$  &  $W$ & $1 + \frac{x}{(1-x)^2}$ & 2 \\
$K \otimes M^*$ & $W^*$ & $\frac{2-x}{(1-x)^2}$ & 1 \\
$W$ & $W$ &  $1 + \frac{x}{(1-x)^2}$ & 2\\
$W \otimes M$ & $F \oplus X_{10}$ & $3 + x +  \frac{x^2}{(1-x)^2}$ & 3\\
$W \otimes M^*$ & $F \oplus F \oplus K$ & $2 + \frac{1}{(1-x)^2}$ & 1\\
$F$ & $F$  & $1$ & 0\\
$F \otimes M$  &   $F^{\oplus 3}$ & $3$ & 0 \\
$F \otimes M^*$ &   $F^{\oplus 3}$ & $3$ & 0 \\
\end{tabular}
\end{center}
\bigskip 

Now we must describe the propagation of these modules in the cohomology
of $({\mathbb Z}_4)^3$.  This follows from the description of $S^*(M)$ in
section~\ref{s:dec-SM}.  It is important to note, however, that we
have doubled all degrees in the symmetric algebra.  For this reason,
we will refer to the invariants as $d_2$, $d_4$, and $d_8$, by abuse
of notation.  We want to compute the Poincar\'e series of the
$E_2$-term, and what we shall do is compute the Poincar\'e series of
the ``tensored with a trivial'' part first, then compute the
Poincar\'e series of the ``tensored with $M$'' part, and then 
compute the Poincar\'e series of the ``tensored with $M^*$'' part.
Finally we combine these pieces.

We make a table for the ``tensored with a trivial part''.
Since there are two trivial modules in $\Lambda^*(M)$, we just consider
the one in degree $0$ and then multiply our Poincar\'e series by
$1+x^3$. In the table below, ``Degree'' means the degree in which the
propagated module first appears in the symmetric algebra, and
``C-Deg.'' means the degree in the $E_2$-term after which we know no
further generators in cohomology appear. 

\bigskip
\begin{center}
\begin{tabular}{c | c | c | c | c}
Module & Propagators & Degree & Poincar\'e Series & C-Deg.\\ \hline
$k$  & $d_8$ &  0 & $\frac{1}{1-x^8} \cdot \frac{1}{(1-x)^2}$ & $0+0=0$ \\
$M$  & $d_8$ & 2 & $\frac{x^2}{1-x^8} \cdot \frac{1}{(1-x)^2}$ & $2+1=3$ \\
$N$  & $d_8$, $d_2$ & 4 &  $\frac{x^4}{(1-x^8)(1-x^2)} \cdot
\frac{1}{1-x}$ & $4+0 = 4$ \\ 
$K$   &  $d_8$, $d_4$ & 4 & $\frac{x^4}{(1-x^8)(1-x^4)} \cdot
\frac{1}{(1-x)^2}$ & $4+1=5$ \\ 
$W$ &  $d_8$, $d_4$ & 6 & $\frac{x^6}{(1-x^8)(1-x^4)} \cdot (1 + 
\frac{x}{(1-x)^2})$ & $6+2=8$ \\ 
$F$ & $d_8$, $d_4$, $d_2$ & 8 & $\frac{x^8}{(1-x^8)(1-x^4)(1-x^2)}$ &
$8+0 = 8$ \\
\end{tabular}
\end{center}
\bigskip 

Thus the total Poincar\'e series for the ``tensored with a trivial''
part is:
\begin{multline*}
(1+x^3) \cdot \biggl[
\frac{1}{1-x^8} \cdot \frac{1}{(1-x)^2} + 
\frac{x^2}{1-x^8} \cdot \frac{1}{(1-x)^2} +
\frac{x^4}{(1-x^8)(1-x^2)} \cdot \frac{1}{1-x} +  \\
\frac{x^4}{(1-x^8)(1-x^4)} \cdot \frac{1}{(1-x)^2} + 
\frac{x^6}{(1-x^8)(1-x^4)} \cdot \bigl(1 + \frac{x}{(1-x)^2}\bigr) + 
\frac{x^8}{(1-x^8)(1-x^4)(1-x^2)} \biggr]
\end{multline*}

Note that the maximum degree of an algebra generator for the
$E_2$-term is $8$.  Although we are tensoring with a trivial 
module in degree
$3$, this is really just multiplying by an element of the $E_2$-term
and so no new algebra generators in degrees $11$ or greater are produced.

Now let's move on to the ``tensored with $M$'' part.

We produce the desired information in the form of a table and a
Poincar\'e series, as above. 

\bigskip
\begin{center}
\begin{tabular}{c | c | c | c | c}
Module & Propagators & Degree & Poincar\'e Series & C-Deg.\\ \hline
$M$  & $d_8$ &  1 & $\frac{x}{1-x^8} \cdot \frac{1}{(1-x)^2}$ & $1+1=2$ \\
$M\otimes M$  & $d_8$ & 3 & $\frac{x^3}{1-x^8} \cdot
\frac{2-2x+x^2}{(1-x)^2}$ & $3+2=5$ \\
$N \otimes M$  & $d_8$, $d_2$ & 5 &  $\frac{x^5}{(1-x^8)(1-x^2)} \cdot
(1+ \frac{1}{1-x})$ & $5+0 = 5$ \\ 
$K \otimes M$   &  $d_8$, $d_4$ & 5 & $\frac{x^5}{(1-x^8)(1-x^4)} \cdot
(1+ \frac{x}{(1-x)^2})$ & $5+2=7$ \\ 
$W \otimes M$ &  $d_8$, $d_4$ & 7 & $\frac{x^7}{(1-x^8)(1-x^4)} \cdot (3+x+ 
\frac{x^2}{(1-x)^2})$ & $7+3=10$ \\ 
$F \otimes M$ & $d_8$, $d_4$, $d_2$ & 9 &
$\frac{3x^9}{(1-x^8)(1-x^4)(1-x^2)}$ & $9+0 = 9$ \\
\end{tabular}
\end{center}
\bigskip 

Thus the total Poincar\'e series for the ``tensored with $M$''
part is:
\begin{multline*}
\frac{x}{1-x^8} \cdot \frac{1}{(1-x)^2} + 
\frac{x^3}{1-x^8} \cdot \frac{2-2x+x^2}{(1-x)^2} +
\frac{x^5}{(1-x^8)(1-x^2)} \cdot (1+ \frac{1}{1-x}) + \\
\frac{x^5}{(1-x^8)(1-x^4)} \cdot (1+ \frac{x}{(1-x)^2}) + \\
\frac{x^7}{(1-x^8)(1-x^4)} \cdot (3+x+ \frac{x^2}{(1-x)^2}) + 
\frac{3x^9}{(1-x^8)(1-x^4)(1-x^2)}
\end{multline*}

Now we do the ``tensored with $M^*$'' part:

\bigskip
\begin{center}
\begin{tabular}{c | c | c | c | c}
Module & Propagators & Degree & Poincar\'e Series & C-Deg.\\ \hline
$k \otimes M^*$  & $d_8$ & 2 & $\frac{x^2}{1-x^8} \cdot
\frac{1}{(1-x)^2}$ & $2+1 = 3$ \\
$M\otimes M^*$  & $d_8$ & 4 & $\frac{x^4}{1-x^8} \cdot ( 1 +
\frac{1}{(1-x)^2})$ & $4+0=5$ \\ 
$N \otimes M^*$  & $d_8$, $d_2$ & 6 &  $\frac{x^6}{(1-x^8)(1-x^2)} \cdot
(1+ \frac{1}{1-x})$ & $6+0 = 6$ \\ 
$K \otimes M^*$   &  $d_8$, $d_4$ & 6 & $\frac{x^6}{(1-x^8)(1-x^4)} \cdot
\frac{2-x}{(1-x)^2}$ & $6+1=7$ \\ 
$W \otimes M^*$ &  $d_8$, $d_4$ & 8 & $\frac{x^8}{(1-x^8)(1-x^4)} \cdot (2+ 
\frac{1}{(1-x)^2})$ & $8+1=9$ \\ 
$F \otimes M^*$ & $d_8$, $d_4$, $d_2$ & 10 &
$\frac{3x^{10}}{(1-x^8)(1-x^4)(1-x^2)}$ & $10+0 = 10$ \\
\end{tabular}
\end{center}
\bigskip 

Thus the total Poincar\'e series for the ``tensored with $M^*$''
part is:
\begin{multline*}
\frac{x^2}{1-x^8} \cdot \frac{1}{(1-x)^2} +
\frac{x^4}{1-x^8} \cdot ( 1 + \frac{1}{(1-x)^2}) + 
\frac{x^6}{(1-x^8)(1-x^2)} \cdot (1+ \frac{1}{1-x}) + \\
\frac{x^6}{(1-x^8)(1-x^4)} \cdot \frac{2-x}{(1-x)^2} + 
\frac{x^8}{(1-x^8)(1-x^4)} \cdot (2+ \frac{1}{(1-x)^2}) + 
\frac{3x^{10}}{(1-x^8)(1-x^4)(1-x^2)}
\end{multline*}

Adding these series up, we obtain the Poincar\'e series for the
$E_2$-term. 
We summarize our computation in the following

\begin{theorem}
The Poincar\'e Series for the $E_2$ term of the Lyndon-Hochschild-Serre spectral
sequence for the cohomology of the group extension $S=Syl_2(HS)=4^3:D_8$
is given by the rational function
$$p(x) ~=~
\frac{(1 +x)^2 (1 - x + x^2) (1 + 2 x^2 - x^5)}
{(1- x)^2(1- x^4)(1- x^8)}.$$
Moreover a complete set of algebra generators for the
$E_2$ term of the spectral sequence 
occur in degrees eight and below.
\end{theorem}

Using Maple, this series can be expanded to yield
\begin{multline*}
1 + 3 x + 7 x^2 + 14 x^3 + 23 x^4 + 34 x^5 + 48 x^6 + 65 x^7 + 84 x^8 +
105 x^9 + 131 x^{10} + \\
163 x^{11} + 198 x^{12} + 236 x^{13} + 280 x^{14} + 330 x^{15}
+383 x^{16} + 439 x^{17} + 503 x^{18} + \cdots
\end{multline*}

Our goal is to analyze the behaviour of this spectral sequence. Our description
provides us with very good control of the generators, especially given that
they arise in low degree. To make effective use of this we require an explicit
computer--assisted calculation of $H^*(S)$ in low degrees. This can indeed be
implemented; details are provided in a subsequent section. For now we use
this information to prove the main result in this paper.

\begin{theorem}\label{specseq}
Let $S=(\mathbb Z/4)^3\rtimes U_3$ denote the Sylow $2$-subgroup of $HS$, 
the Higman--Sims group. The Lyndon--Hochschild--Serre spectral sequence for
this semidirect product collapses at $E_2$ and hence the polynomial
$p(x)$ is the Poincar\'e series for $H^*(S,\mathbb F_2)$.
\end{theorem}
\begin{proof}
A computer calculation using MAGMA for $H^*(S,\mathbb F_2)$
shows that through degree 10 the coefficients of the polynomial
$p(x)$ agree with the ranks of the cohomology.
Now according to the preceding
tables, all algebra generators occur by this degree. Hence using the multiplicative
structure of the spectral sequence we infer that it must collapse at $E_2$, i.e.
$E_2=E_{\infty}$.
\end{proof}

\section{Computer Calculations of the Cohomology of $S=Syl_2(HS)$}

The cohomology $H^*(S,\mathbb F_2)$ was calculated directly 
through degree 10 using a computer. The calculation included a 
determination of the Betti numbers, a minimal set of ring generators, 
and a complete collection of relations among the generators in the
first ten degrees. We also obtained the images of the generators under
various restrictions. Using the latter information and the analysis
of the last section we were able to construct the cohomology ring
with all generators and relations using computer technology. 
We give a summary of the calculations in this section. 

The calculation was begun by first obtaining a minimal projective
resolution 
$$
\ldots \longrightarrow P_2 \longrightarrow P_1  \longrightarrow P_0
 \longrightarrow k
$$
of the trivial module $k =  \mathbb F_2$. The process is 
mostly linear algebra. The free module over $S$ is generated as
a collection of matrices representing the actions of the generators
of the group. At each stage in the construction we get a minimal
set of generators for the kernel 
of the previous boundary map, create the free module with
exactly that number of generators, make the matrix for the boundary
map from the free module to the kernel of the previous boundary map
-- this is the new boundary map -- and find its kernel. The programs
are conservative with both time and memory in that they save only
the minimal amount of information necessary to reconstruct the 
boundary homomorphisms and create module structure only when it is 
necessary. The method is described in detail in \cite{CGS},
\cite{Ctest}.

Because the resolution is minimal, any nonzero homomorphism 
$\gamma^{\prime}: P_n \longrightarrow k$ is a cocycle representing a nonzero 
cohomology class $\gamma$.  Any such cocycle can be lifted to a chain map
$\hat\gamma: P_* \longrightarrow P_*$ of degree $n$ which also represents
$\gamma$. Again this is a exercise in linear algebra on 
the computer and details of the implementation can be found in the
references given above. The point of this operation is that the cup 
product of two cohomology elements is the class of the composition of 
the representing chain maps. Thus we can compute the product structure
of the cohomology ring. Note that the program computes the chain maps
for only a minimal set of generators for the cohomology. That is, 
at each stage it constructs the subspace of $\Hom_{kG}(P_n,k)$ that
can be obtained by compositions with chain maps of lower degrees.

The same implementations have been used to compute the cohomology
rings of all but a few of the groups of order dividing $64$.\footnote{See
the second author's web page, 
http://www.math.uga.edu/\symbol{126}jfc/groups2/cohomology2.html
for the results and some further
discussion of the methods.} The programs are written in the MAGMA 
language and run in the MAGMA computer algebra system \cite{BCann}.
The computations of the cohomology of $S$ were all run on an
SUN ULTRA2200 (named the sloth) that has approximately 1Gb. of RAM
and 14Gb. of hard disk. The computation of the projective resolution
out to degree $10$ for the trivial module $k$ of $S$ took 
slightly under $31$ hours. The computation of the chain maps of the
minimal generators of the cohomology ring took more than $55$ hours. 
The attempt to compute an $11^{th}$ step in the projective resolution
failed for lack of memory. Some of the results of this calculation
are given in the following. Note that this is precisely the information
that is needed to complete the proof of Theorem \ref{specseq}.

\begin{proposition}
\label{betti}
The Betti numbers for the cohomology $H^*(S,k)$ through 
degree $10$ are 
$$ 
1, 3, 7, 14, 23, 34, 48, 65, 84, 105, 131,
$$
and the degrees of the minimal generators of the cohomology 
ring through degree $10$ are 
$$
1, 1, 1, 2, 2, 2, 2, 3, 3, 3, 3, 4, 4, 5, 6, 7, 8.
$$
\end{proposition}

We also computed the restrictions of the cohomology to the elementary
abelian subgroups and to the centralizers of the elementary 
abelian subgroups. The restriction maps are obtained by first
converting the projective resolution of $k$ for the group to a projective 
resolution of the subgroup and then constructing the chain map,
lifing the identity on $k$, from
a minimal projective resolution of $k$ for the subgroup to the converted
projective resolution for the group. The cocycles for the cohomology
elements of the group are then pulled back along the chain map. The 
process is described in more detail in \cite{CMM}.

Of particular interest was the restrictions to the centralizers of the 
elementary abelian groups of order $4$. This calculation was made on the
assumption (hope) that the depth of the cohomology ring $H^*(S,k)$ 
is at least two and hence that the cohomology of $S$ is detected
on these centralizers. An easy calculation shows the following.

\begin{lemma}
If $E$ is an elementary abelian subgroup of order four in $S$ 
then the centralizer of $E$ is contained in the subgroup $K_{\beta}$
(defined in Section ~2) or in one of the three subgroups
$$
M = \langle v_1,v_2,v_3, ts, s^2 \rangle,  
$$
$$
D_1 = \langle (v_1v_2)^2, v_1v_3,t,s \rangle,  
$$
$$
D_2 = \langle (v_1v_2)^2, v_1v_3, v_2^{-1}t,ts,v_1v_2^{-1}s^2 \rangle.
$$
Here $M$ is the centralizer of $(v_1v_2)^2$, $D_1$ is the centralizer
of $s^2$ and $D_2$ is the centalizer of $v_1v_2^{-1}s^2$. The orders
of $M$ is  $256$ while the orders of $D_1$ and $D_2$ are $64$.
\end{lemma}

The computer was able to calculate the cohomology ring of all four
of the groups in the Lemma. In the case of $K_{\beta}$ it was known
that the cohomology is detected on the centralizers of the maximal
elementary abelian subgroups and this aided the computation. In all
three of the other cases the minimal set of generators for the 
cohomology ring contained no element of degree greater than $4$ and
the computation was possible even though all three of the groups have
$2$-rank $4$ and the Betti numbers of the cohomology grow rapidly.
In fact, the cohomology rings of $D_1$ and $D_2$ are included in the
calculations of the second author (see the web page). These groups 
have Hall-Senior numbers 170 and 110 respectively. 

We were able to compute the restriction maps and then get the 
kernels of the restriction to each of the subgroups. Then the 
intersection of the kernel of the restriction maps was computed.
Actually the problem of getting the intersection of the kernels
seemed to be too difficult for the computer to attack directly
using Gr\"obner basis machinery.
Instead we employed an indirect method of turning the restriction
maps into linear transformations on the spaces of monomials in
each degree and computing the intersection of the kernels out to
degree $17$. Thus we had a complete set of generators for the intersection
of the kernels out through degree $17$. 

There is one uncertainty that should be noted. The cohomology of $M$ 
was only calculated out to degree $10$ and that is not far enough 
to pass our test for completeness of the calculation \cite{Ctest}. 
Nonetheless, the generators of the cohomology of $S$ are in degrees
at most $8$ and so their restrictions to $M$ are expressed as 
polynomials in the calculated generators of the cohomology of $M$. 
Any polynomial in those generators which is computed to be in the 
kernel of restriction to $M$ is, in fact, in that kernel. Thus we get
the following result, which is used in the next section.

\begin{proposition}
The ideal ${\mathcal I}$ of Theorem \ref{co-sy2} is generated
by the elements of degree at most $17$ 
that are in the intersection of the kernels 
of the restriction maps to $K_{\beta}$, $M$, $D_1$ and $D_2$. 
Hence the ideal is in the intersection of the kernels of the 
restrictions to the centralizers of the elementary abelian subgroups
of order four.
\end{proposition}

In the next section we argue that this information together with some
extra verification is sufficient.

\section{Generators and Relations for the Cohomology of $S$}

Our purpose in this section is to give a proof of the following.

\begin{theorem}
\label{co-sy2}
The cohomology of 
$S$ has the form
$$H^*(S,{\mathbb F}_2) \cong 
{\mathbb F}_2[z,y,x,w,v,u,t,s,r,q,p,n,m,k,j,i,h]/{\mathcal I}$$
where the degree of $z, y, x$ is $1$, the degree of $w, v, u, t$ 
is $2$, the degree of $s, r, q, p$ is $3$, the degree of $n, m$
is $4$ and the degrees of $k,j,i,h$ are $5, 6, 7$ and $8$ respectively.
The ideal ${\mathcal I}$ is generated by the relations 
in the appendix.
\end{theorem}

Let $R = {\mathbb F}_2[z,y,x,w,v,u,t,s,r,q,p,n,m,k,j,i,h]$ be the 
polynomial ring in $17$ variables in weighted degrees as in the 
theorem. For $n \geq 0$ let $R^n$ be the space of homogenous 
polynomials of degree $n$. From Proposition \ref{betti} we know
that there is a natural homomorphism $\phi: R \longrightarrow
H^*(S,{\mathbb F}_2)$ taking each variable to the computed
cohomology element. Let ${\mathcal J}$ be the kernel
of $\phi$. Our task is
to establish that ${\mathcal I} \subseteq {\mathcal J}$.

Recall that the
following is a complete list of
representatives for the conjugacy classes of centralizers of
maximal elementary abelian subgroups in $S$: 
\begin{enumerate}
\item $CI~=~4^3$, a characteristic subgroup of $S$
\item $CII~=~2^2\times 4$,
\item $CIII~=~2^2\times 4$,
\item $CIV~=~2^2\times 4$,
\item $CV~=~2^2\times 4$,
\item $CVI~=~2^2\times 4$, 
\item $2^4_A$, 
\item $2^4_B$
\item $2^4_C$
\end{enumerate}

The following basic detection result involving these groups, these 
groups was proved by computer calculation. 

\begin{lemma}
Any homogeneous polynomial of $R^n$  in degree 
$n \leq 20$ is in ${\mathcal I}$ if and only if it is in 
the kernel of restriction to all of the centralizers 
listed above.
\end{lemma}
\begin{proof}
First it was shown that the elements of ${\mathcal I}$ are in the kernels 
of restriction by computing the restriction maps on the elements. Then
for each degree $n \leq 20$ the restriction maps $R^n/{\mathcal I}^n
\longrightarrow H^*(C,{\mathbb F}_2)$ for $C$ in the list was written
as a linear transformation and the intersection of the null spaces 
was computed to be zero
\end{proof}

\begin{proof} [Proof of Theorem \ref{co-sy2}] 
 From the lemma we see
that ${\mathcal J}^n \subseteq {\mathcal I}^n$ for $n \leq 20$.
But by the Poincar\'e series (Theorem \ref{specseq})
we have that $R^n/{\mathcal I}^n$
has the same dimension as $R^n/{\mathcal J}^n \cong H^*(S,{\mathbb F}_2)$.
On the other hand we know that ${\mathcal I}$ is generated by elements
of degree at most $14$. So ${\mathcal I} \subset {\mathcal J}$, since
all generators are in ${\mathcal J}$. So $\phi$ induces a 
surjective homomorphism
$R/{\mathcal I} \longrightarrow H^*(S,{\mathbb F}_2)$. But again 
because the Poincar\'e series are the same this is an isomorphism.
\end{proof}
                         
\section{Appendix I: Relations in the cohomology of $S$}

The following is a list of the relations in the cohomology ring
$H^*(S,{\mathbb F}_2)$. These elements generate the ideal ${\mathcal I}$
in Theorem \ref{co-sy2}.
Note that there are ~79 relations generating ${\mathcal I}$ and 
that the relations are minimal in the sense that no collection of 
fewer than ~79 elements will generate the ideal. This set is not a 
Gr\"obner basis for the ideal. The computed Gr\"obner basis using
the grevlex ordering on the variables consisted of 884 elements 
in degrees up to $73$. The elements in the list are ordered by 
degree. The largest degree of any element in the list is $14$.
\vskip.2in
$\phantom\qquad    
yx,
~~zx,
~~zy,
~~xv,
~~xw,
~~yw + yv + yu,
~~zt,
~~zw,
~~vt,
~~xs,$\newline

$\phantom\qquad    
yq + xq + w^2 + wu,
~~yq + wv,
~~yr + xr + xq + xp + ut,
~~yr + wt,
~~ys + yq,$\newline

$\phantom\qquad    
y^2u + yq + xq + wu,
~~zp + yq + xq + v^2 + vu,
~~zr,
~~zs + zq + xq,$\newline

$\phantom\qquad    
tq,
~~ts,
~~vr,
~~vs + vq,
~~ws + wq,
~~x^2r + xn + wp + vp + up,$\newline

$\phantom\qquad    
x^2r + xn + wr + ur,
~~yu^2 + yn + wp,
~~yvu + ws,$\newline

$\phantom\qquad    
y^2r + wr,
~~z^2p + zvu + zn + ws + vs + us,$\newline

$\phantom\qquad    
z^3u + z^2p + zvu + zu^2 + zn + zm + x^2r + xu^2 + xn + xm,$\newline

$\phantom\qquad    
rq,
~~sp + qp,
~~sq + q^2,
~~sr,
~~s^2 + sq,
~~xuq + vu^2 + vn + sp,$\newline

$\phantom\qquad    ytr + x^3p + x^2n + xup + xtp + r^2,
~~ytr + x^3p + xuq + xtr + xk + r^2,$\newline

$\phantom\qquad    yur + x^3p + xuq + xk + tn + rp,
~~y^2n + yvp + yuq + wn,$\newline

$\phantom\qquad    z^3p + z^2n + zvq + zuq + y^2m + ytr + ytp + yk + x^3p + xk + 
        vn + vm + tn + tm + p^2,$\newline

$\phantom\qquad    z^3p + zvq + zup + yuq + xuq + vu^2 + sq,
~~z^3q + z^3p + zuq + zk,
~~sn + qn,$\newline

$\phantom\qquad    ytn + xtm + xp^2 + rn,
~~yvn + yqp + vuq + vup + qn,$\newline

$\phantom\qquad    y^5v + y^4q + y^4p + y^3n + y^3m + y^2vq + y^2tp + y^2k + yum + 
        ytn + yr^2 + yrp + yj + wk,$\newline

$\phantom\qquad    zvm + zq^2 + zj + yvm + yum + yr^2 + yrp + x^3u^2 + x^3m + 
        vup + vk + u^2q + u^2p + uk + sm + rn + rm + pn,$\newline

$\phantom\qquad    z^2vq + zvm + zq^2 + zqp + zj + yvn + yqp + vuq + u^2s,$\newline

$\phantom\qquad    z^3u^2 + zu^3 + zum + yvn + yqp + x^3u^2 + x^3n + xu^3 + xum + 
        xtn + xj + vuq + u^2s,$\newline

$\phantom\qquad    z^5u + z^3m + zu^3 + zum + x^3n + x^3m + xu^3 + xum + xtn + xj 
        + wk + vk + u^2s + u^2q + uk,$\newline

$\phantom\qquad    y^4m + y^3k + y^2vn + y^2vm + y^2p^2 + yt^2p + ytk + yrm + 
        x^4u^2 + x^3k + x^2un + x^2tn + x^2j + xu^2p + xtk + xrm + 
        xpn + tr^2 + tj + rk,$\newline

$\phantom\qquad    y^5p + y^3k + y^2vn + y^2vm + yrm + ypm + yi + x^4u^2 + x^3k + 
        x^2un + x^2tn + x^2j + xt^2p + xqm + wun + wum + wj + 
        tr^2 + sk + qk,$\newline

$\phantom\qquad    zpn + ypn + xpn + wun + vun + vp^2 + uqp + n^2,$\newline

$\phantom\qquad    z^2q^2 + zuk + zqn + zqm + y^4m + y^3vq + y^3tp + y^3k + 
        y^2vm + y^2tm + y^2p^2 + yqm + ypn + x^2un + xuk + xrm + 
        wun + wj + vum + vp^2 + vj + up^2 + qk,$\newline

$\phantom\qquad    z^3k + z^2vm + z^2q^2 + zi + y^3tp + y^2tm + yuk + yt^2p + 
        ytk + yqm + x^3k + x^2u^3 + x^2un + x^2tn + xt^2p + xi + 
        wum + vun + vum + urp + up^2 + tr^2 + tj + sk + rk + qk,$\newline

$\phantom\qquad    z^3k + z^2vm + zqn + zpn + zi + y^5q + y^4n + y^2q^2 + yqm + 
        x^2u^3 + xu^2p + xuk + xrm + wj + vq^2 + vp^2 + u^2n + uq^2 
        + uj + sk,$\newline

$\phantom\qquad    z^4u^2 + z^2u^3 + z^2un + z^2um + z^2q^2 + zqn + y^5q + y^4n + 
        y^2q^2 + yqm + wun + vun + vq^2 + vp^2 + vj + uq^2 + sk,$\newline

$\phantom\qquad    z^5p + z^3k + z^2un + z^2j + zqn + zi + y^3vq + y^2vn + yuk + 
        yrm + ypn + x^2u^3 + xuk + xqm + wum + vp^2 + u^2n + urp 
        + uj + sk,$\newline

$\phantom\qquad    y^5n + y^3q^2 + y^2vk + y^2qm + yvqp + yup^2 + yuj + yt^2n + 
        yrk + yqk + ynm + wuk + wqn + wpn + wi + trm + sj + 
        r^2p + rp^2 + qj,$\newline

$\phantom\qquad    z^2pn + zuq^2 + zuj + zn^2 + znm + y^6p + y^5m + y^3vm + 
        y^3tm + y^2qn + y^2qm + y^2i + yvqp + yvj + yt^2n + 
        yt^2m + ytj + yrk + ypk + ynm + wqn + wpm + wi + vpm + 
        u^3q + u^2k + urn + trn + q^3 + pj + nk,$\newline

$\phantom\qquad    z^2qn + zuj + zn^2 + znm + yup^2 + yt^2n + yrk + ynm + xpk 
        + wpn + wpm + vuk + vqm + u^3q + u^2k + upn + trm + sj + 
        rp^2 + rj + qj,$\newline

$\phantom\qquad    z^3un + z^2qn + z^2pn + zvj + zuq^2 + zuj + znm + sj + qj,$\newline

$\phantom\qquad    y^2tk + y^2qm + yup^2 + yt^2m + ytp^2 + ytj + yqk + ynm + 
        x^3tn + xrk + xn^2 + wuk + wqn + wpm + t^3r + trn + tpm 
        + ti + sj + r^2p + rp^2 + qj,$\newline

$\phantom\qquad    z^3vn + z^3un + z^2vk + z^2pn + zvj + y^5n + y^4vq + y^3vm + 
        y^2qn + y^2pn + yvq^2 + yvqp + yup^2 + yn^2 + xnm + wuk 
        + vqm + vpm + vi + u^3q + u^3p + u^2k + upn + q^3,$\newline

$\phantom\qquad    z^3vn + z^2qn + z^2pn + zvj + y^5n + y^2vk + y^2tk + y^2qn + 
        yvqp + yup^2 + yt^2n + yt^2m + ytp^2 + ytj + yrk + yn^2 
        + x^3tn + wpn + wpm + wi + vuk + vqm + urn + t^3r + trm 
        + tpm + ti + sj + r^2p,$\newline

$\phantom\qquad    y^3t^2p + y^2t^2m + yt^3p + yt^2k + yr^2p + yrp^2 + x^2u^2n + 
        x^2t^2n + xrj + t^3n + t^2rp + t^2j + trk + r^2n + r^2m + 
        rpm + ri,$\newline

$\phantom\qquad    z^2uq^2 + z^2qk + zpj + y^5vq + y^4j + y^3vk + y^3t^2p + 
        y^3qn + y^3qm + y^3pm + y^3i + y^2vq^2 + y^2vj + y^2t^2m + 
        y^2tp^2 + y^2tj + y^2pk + y^2h + yvpn + yupn + yt^3p + 
        yt^2k + yr^2p + yqp^2 + yqj + x^2t^2n + xu^3p + xtpn + 
        xrj + xpj + wqk + wpk + wnm + vup^2 + upk + t^3n + 
        t^2p^2 + t^2j + trk + r^2n + r^2m + rpn + rpm + q^2n + p^2n 
        + pi + k^2,$\newline

$\phantom\qquad    zuqn + y^7q + y^4q^2 + y^3vk + y^3qm + y^2vq^2 + y^2vj + 
        yvqm + yvpn + yvi + yupn + yqj + x^2u^2n + x^2nm + 
        xu^3p + xpj + wqk + wpk + vuj + upk + un^2 + tn^2 + r^2n 
        + q^2n + qpn + qpm + nj,$\newline

$\phantom\qquad    z^3uk + zuqm + zupm + zmk + y^3vk + y^3qn + y^3pm + 
        y^2tp^2 + y^2qk + yvpn + yvi + yupn + yupm + yr^2p + 
        yq^2p + yqj + ypj + x^2u^2n + xpj + wpk + vup^2 + upk + 
        t^2rp + rpn + qpm + p^2n + p^2m + pi,$\newline

$\phantom\qquad    z^3uk + z^2u^2n + z^2uj + z^2qk + zupm + zui + zpj + y^3qn 
        + y^3qm + y^3pn + y^2vq^2 + y^2nm + yq^2p + x^3uk + x^2u^4 
        + x^2nm + xu^3p + xui + xrj + wpk + wnm + vuj + usk + 
        upk + un^2 + tn^2 + si + r^2n + qpn + nj,$\newline

$\phantom\qquad    z^4vm + z^4j + z^3vk + zuqn + zuqm + zupm + zpj + y^7q + 
        y^4q^2 + y^3vk + y^3qn + y^3pn + y^2vj + y^2nm + yvpn + 
        yvi + yq^2p + yqj + x^2u^2n + x^2nm + vup^2 + u^3n + usk 
        + un^2 + unm + qpn + qi,$\newline

$\phantom\qquad    z^2u^2k + z^2uqm + z^2q^3 + z^2qj + zvqk + zu^3n + zu^2j + 
        zupk + zun^2 + y^6vq + y^4qm + y^4pn + y^3tp^2 + y^2qj + 
        yvpk + yvn^2 + yvnm + yupk + yq^2n + yqpm + yqi + 
        x^2u^2k + x^2nk + xunm + xtpk + xtnm + xri + xnj + 
        wupm + wpj + wnk + vq^2p + vpj + urp^2 + unk + t^2rn + 
        t^2pm + trp^2 + tp^3 + tpj + rpk + qpk + qnm + pnm + ni,$\newline

$\phantom\qquad    z^2vqm + z^2u^2k + z^2uqm + z^2q^3 + zvqk + zu^3n + zu^2j + 
        zun^2 + y^9v + y^6vq + y^6k + y^5vm + y^5tm + y^5q^2 + 
        y^4qn + y^4qm + y^3tp^2 + y^3tj + y^3pk + y^3nm + y^3m^2 + 
        y^3h + y^2vi + y^2qj + y^2nk + yvpk + yvnm + yvm^2 + 
        yvh + yupk + yuh + ytm^2 + x^2mk + xunm + xtnm + xqi 
        + wupm + wpj + wnk + wmk + vmk + urj + uqp^2 + upj + 
        unk + t^2pm + tp^3 + r^2k + rn^2 + qpk + qn^2 + qnm + ni + 
        kj,$\newline

$\phantom\qquad    z^2q^2n + zvqj + zvmk + zupj + zqn^2 + zni + y^8m + y^7vq 
        + y^7tp + y^7k + y^6tm + y^5vk + y^5qn + y^4pk + y^4nm + 
        y^3p^3 + y^2vnm + y^2vh + yvq^3 + yvqj + yvnk + yvmk + 
        yup^3 + yumk + ytrj + ytnk + ytmk + yrpk + yqm^2 + 
        yp^2k + ypnm + yni + ykj + x^2mj + xqm^2 + xkj + wum^2 
        + wq^2m + wpi + vqpn + vp^2m + u^2qk + u^2pk + u^2n^2 + 
        urpm + uri + uq^2m + up^2n + t^2nm + t^2m^2 + tr^2n + 
        tp^2m + rpj + q^4 + q^3p + q^2j + qp^3 + p^2j + pnk + j^2,$\newline

$\phantom\qquad    z^5pm + z^2q^2n + z^2qi + zvqj + zupj + zpm^2 + zkj + y^8m 
        + y^7vq + y^7k + y^6vm + y^5tk + y^5pm + y^4pk + y^3vqn + 
        y^3t^3p + y^3ti + y^2vnm + y^2vh + y^2th + y^2q^2m + 
        y^2k^2 + yvqj + yvnk + yumk + yt^4p + yt^3k + ytnk + 
        ytmk + yrpk + yrn^2 + yq^2k + yqn^2 + yp^2k + ypnm + 
        yni + ykj + x^8u^2 + x^8m + x^7k + x^6u^3 + x^5tk + x^5i + 
        x^4tj + x^2t^2j + x^2mj + xt^2pn + xni + wun^2 + wum^2 + 
        wq^2m + wmj + vum^2 + vqpn + vnj + vmj + u^2qk + uri + 
        uqi + up^2n + upi + uk^2 + t^4n + t^3rp + t^3p^2 + t^3j + 
        t^2nm + t^2m^2 + trpm + tk^2 + r^3p + r^2p^2 + rpj + q^4 + 
        q^3p + q^2j + qp^3 + p^2j + pmk + ki + j^2,$\newline

$\phantom\qquad    z^{11}u + z^9u^2 + z^9m + z^7vn + z^6vk + z^6pm + z^3vnm + 
        z^3un^2 + z^3nj + z^2unk + z^2pm^2 + z^2mi + z^2kj + 
        zumj + y^7p^2 + y^6vk + y^4pj + y^4nk + y^3qi + y^3pi + 
        y^2vnk + y^2tmk + y^2q^2k + y^2qnm + y^2qh + y^2p^2k + 
        yvq^2m + yvqpm + yvqi + yvk^2 + yup^2n + yup^2m + 
        yumj + yuk^2 + ytri + yr^2p^2 + yr^2j + yq^4 + yq^2j + 
        yqpj + yn^2m + x^6uk + x^5u^2m + x^4ui + x^4mk + x^3nj + 
        xu^4n + xu^2nm + xupi + xumj + xtnj + xpnk + wunk + 
        wq^2k + wmi + vupj + vumk + vqnm + u^5q + u^4k + u^2p^3 
        + u^2nk + usm^2 + uq^2k + uqn^2 + t^3rn + trnm + tpnm + 
        tpm^2 + q^3n + qp^2n + qnj + p^3m + pnj + nmk + ji,$\newline

$\phantom\qquad    z^{11}p + z^5vi + z^5u^2k + z^4pi + z^3umk + z^2qmk + z^2j^2 + 
        zvqn^2 + zuqnm + zumi + zukj + zq^3n + zm^2k + y^{10}n +
        y^9tp + y^8q^2 + y^8j + y^7tk + y^7qn + y^7i + y^6vj + 
        y^6tj + y^6nm + y^6h + y^5vi + y^5ti + y^4vn^2 + y^4vm^2 + 
        y^4mj + y^4k^2 + y^3vqj + y^3tmk + y^3p^2k + y^3mi + 
        y^2vqi + y^2vnj + y^2t^2m^2 + y^2q^2j + y^2n^2m + y^2mh + 
        y^2ki + yvq^2k + yvqm^2 + yvqh + yvpnm + yvmi + 
        yupn^2 + yuph + yuni + yt^3pm + yt^2pj + ytp^2k + 
        ytmi + yr^2pm + yq^3n + yp^2i + yn^2k + x^8u^3 + x^7uk + 
        x^6u^2m + x^6uj + x^4um^2 + x^3mi + x^3kj + x^2u^4n + 
        x^2umj + x^2j^2 + xu^2mk + xrk^2 + xji + wunj + wqnk + 
        wnm^2 + wnh + vumj + vq^2p^2 + vqmk + vn^2m + vj^2 + 
        u^2qi + u^2p^2m + u^2mj + urpj + urmk + uqnk + upmk + 
        un^3 + unm^2 + uki + t^3nm + t^2r^2m + t^2rpm + t^2ri + 
        t^2p^2n + tr^2j + trp^3 + trnk + tp^4 + skj + r^3k + 
        r^2pk + rpn^2 + rni + qmi + pni + mk^2 + i^2$\newline

\section{Appendix II: Steenrod Operations}
In this appendix we describe the Steenrod operations on the cohomology generators
listed previously. As before this was done on a computer using MAGMA. Note that
the program used puts all of the polynomials in ```normal form'' relative to
the Groebner basis of the ideal of relations.
\vskip .2 in

$\phantom\qquad
Sq^1w=yv, 
~~Sq^1v=zv+yv,
~~Sq^1u=zu+xu,
~~Sq^1t=yt+xt$\newline

$\phantom\qquad
Sq^1s=vu+yq+xq,~~~Sq^2s=z^2q+us+zn$\newline

$\phantom\qquad
Sq^1r=xr+xp,~~~Sq^2r=ur+tr+xn$\newline

$\phantom\qquad
Sq^1q=vu+yq+xq,~~~Sq^2q=xu^2+z^2q+uq+wp+vp+up+zn+xm$\newline

$\phantom\qquad
Sq^1p=zp,~~~Sq^2p=z^2q+z^2p+y^2p+tr+wq+uq+vp+ym+k$\newline

$\phantom\qquad
Sq^1n=x^2p+wp+vp,~~~Sq^2n=y^4v+y^3q+y^3p+zvq+yvq+xuq+yvp+
			zup+yup+ytp+z^2n+y^2n+x^2n+y^2m+
			x^2m+q^2+rp+qp+tn+vm+yk+xk+j$\newline

$\phantom\qquad
Sq^1m=zu^2+xu^2+y^2p+
	us+vq+up+zm+xm$\newline

$\phantom\qquad
Sq^2m=  y^4v+y^3q+y^3p+u^3+zvq+xuq+
			   yvp+xup+xtp+z^2m+y^2m+r^2+q^2+
			   rp+vn+un+wm+vm+um+tm+xk+j$\newline

$\phantom\qquad
Sq^1k=y^3p+ytr+zvq+yvq+yvp+yup+xtp+z^2n+y^2m+qp+p^2+vn$\newline

$\phantom\qquad
Sq^2k=y^4p+xu^3+x^3m+u^2s+u^2p+yrp+yp^2+zun+yun+zvm+yvm+
xum+xtm+y^2k+x^2k+rn+qn+pn+sm+rm+vk+xj$\newline
         
$\phantom\qquad
Sq^4k=xu^4+y^4k+yvq^2+xt^2n+z^2pn+xu^2m+yt^2m+y^2vk+y^3j+r^2p+qp^2+p^3+
urn+uqn+tpn+zn^2+uqm+wpm+upm+znm+xnm+ym^2+wuk+vuk+u^2k+t^2k+yrk+
xrk+zqk+yqk+xpk+xuj+ytj+xtj+z^2i+x^2i+mk+pj+vi+yh
$\newline   

$\phantom\qquad
Sq^1j=y^4q+y^2vq+y^3n+y^3m+x^3m+zq^2+vup+zqp+yvn+ytn+zvm+ytm+y^2k+qn+pn+vk+zj+yj+xj
$\newline

$\phantom\qquad
Sq^2j=y^2p^2+x^2tn+z^2vm+ur^2+vq^2+tp^2+t^2n+zqn+vum+t^2m+zqm+xqm+zpm+
zvk+yvk+zuk+xtk+z^2j+n^2+sk+vj
$\newline

$\phantom\qquad
Sq^4j=y^2t^2m + y^3pm + z^4j + x^4j + yq^3 + t^2rp + yqp^2 + yp^3 + 
    yvqn + yupn + y^2n^2 + yvqm + yupm + z^2nm + xu^2k + 
    z^2qk + z^2vj + y^2tj + x^2tj + y^3i + qpn + un^2 + q^2m + 
    qpm + wnm + tnm + wm^2 + vm^2 + tm^2 + urk + trk + wqk + 
    uqk + wpk + upk + znk + xnk + ymk + vuj + yrj + zqj + 
    zpj + xpj + xti + y^2h + nj + mj + qi + pi + wh$\newline

$\phantom\qquad
Sq^1i=y^2q^2+y^2p^2+z^2vm+y^2tm+ur^2+tr^2+wq^2+urp+trp+
    vqp+vp^2+tp^2+vun+xrn+yqn+zpn+xpn+vum+
    t^2m+yqm+xqm+zpm+ypm+zuk+xuk+z^2j+n^2+sk+
    rk+wj+vj+tj+yi$\newline

$\phantom\qquad
Sq^2i=y^3tm+yvq^2+xt^2n+y^2pn+y^2tk+q^2p+rp^2+qp^2+p^3+
    vqn+uqn+zn^2+xn^2+urm+trm+wqm+uqm+wpm+
    vpm+upm+tpm+ynm+xnm+ym^2+xm^2+vuk+xrk+
    zqk+xpk+zvj+xuj+z^2i+x^2i+nk+mk+sj+rj+qj 
   +vi+ti$\newline

$\phantom\qquad
Sq^4i=y^4tk+xu^3m+yt^3m+y^3vj+y^3tj+z^4i+y^4i+x^4i+vq^3 
   +vupn+yqpn+yvn^2+yun^2+t^2pm+yqpm+xunm+
    yum^2+xum^2+ytm^2+u^3k+xupk+z^2nk+y^2nk+z^2mk 
   +y^2mk+yt^2j+xt^2j+z^2qj+y^2qj+z^2pj+y^2vi+
    y^2ti+rn^2+pn^2+rnm+rm^2+qm^2+r^2k+vnk+tnk+
    vmk+umk+tmk+yk^2+urj+wqj+vqj+uqj+wpj+
    ynj+xnj+zmj+ymj+xmj+t^2i+yqi+ypi+xpi+
    yvh+kj+ni+mi+sh+qh$\newline

$\phantom\qquad
Sq^1h=z^7v+x^7u+y^2qn+z^2pn+y^2pn+yt^2m+y^2pm+y^2vk+
    y^2tk+x^2tk+y^3j+rp^2+p^3+vpn+upn+tpn+xn^2+
    usm+trm+wqm+vqm+wpm+vpm+xnm+wuk+vuk+
    yrk+zpk+ypk+xpk+zvj+yvj+zuj+yuj+xuj+
    x^2i+nk+pj$\newline
         
$\phantom\qquad
Sq^2h=z^8v+x^8u+u^5+y^3pm+y^3tk+t^2rp+yrp^2+yqp^2+yupn 
   +y^2n^2+u^3m+yvqm+ytpm+z^2nm+y^2nm+x^2nm+
    z^2m^2+y^2m^2+x^2m^2+xu^2k+xt^2k+z^2qk+y^2vj+
    x^2uj+y^2tj+x^2tj+x^3i+r^2n+q^2n+p^2n+wn^2+vn^2 
   +un^2+q^2m+rpm+qpm+wnm+vnm+trk+wqk+uqk+
    wpk+upk+ynk+xnk+zmk+vuj+zqj+ypj+yvi+
    xui+yti+xti+y^2h+pi$\newline

$\phantom\qquad
Sq^4h=z^6j+x^6j+u^6+y^2t^3m+y^3t^2k+x^5i+t^3p^2+t^4n+t^4m+
    z^2um^2+xt^3k+y^3nk+y^3pj+y^3vi+x^3ti+z^4h+y^4h 
   +x^4h+q^4+p^4+vq^2n+tp^2n+u^2n^2+yrn^2+yqn^2+
    ypn^2+tr^2m+vq^2m+vqpm+vp^2m+tp^2m+vunm+
    u^2nm+t^2nm+yrnm+ypnm+wum^2+vum^2+u^2m^2+
    yqm^2+xqm^2+ypm^2+zq^2k+yrpk+yqpk+zvmk+
    zumk+yumk+u^3j+t^3j+yvpj+xupj+ytpj+xtpj+
    y^2nj+x^2nj+z^2mj+xu^2i+xt^2i+z^2qi+y^2qi+
    y^2pi+x^2pi+y^2vh+z^2uh+nm^2+m^3+pmk+wk^2+vk^2 
   +r^2j+q^2j+rpj+vnj+unj+tnj+wmj+tmj+zkj+
    ykj+xkj+usi+uri+vqi+uqi+wpi+upi+zni+
    yni+ymi+vuh+t^2h+yrh+xrh+zph+ki+nh+mh$\newline

\end{document}